\documentclass[11pt,a4paper]{article}

\usepackage{cite}
\usepackage[cmex10]{amsmath}
\usepackage{amsmath}
\usepackage{amssymb}
\usepackage{booktabs}
\usepackage{color}
\usepackage{times}
\usepackage{paralist}
\usepackage{subfigure}
\usepackage{float}
\usepackage{epstopdf}
\usepackage{authblk}
\usepackage{graphicx}

\let\mathnumsetfont\mathbb
\newcommand\Prob{\mathnumsetfont P}
\newcommand\Mean{\mathnumsetfont E}
\newcommand\Unif{\mathnumsetfont U}
\newcommand\Normal{\mathnumsetfont N}

\usepackage[utf8]{inputenc}
\usepackage[cmex10]{amsmath}
\usepackage{amssymb}
\usepackage{mathtools}
\usepackage{subfigure}
\usepackage{float}
\usepackage{booktabs}
\usepackage{times}
\usepackage{paralist}
\usepackage{cite}
\usepackage{algorithm}
\usepackage{algpseudocode}
\usepackage{epstopdf}
\usepackage{soul} % allows for  \hl{}
\usepackage{color}
\usepackage{eurosym}
\usepackage[printonlyused]{acronym}
\graphicspath{
    {FiguresBldsExp/}
}
\DeclareGraphicsExtensions{.pdf,.png,.jpg,.jpeg,.eps}

\usepackage{amsthm}

\newtheorem{proposition}{Proposition}

\newtheorem{remark}{Remark}

\definecolor{Red}{rgb}{1,0,0}
\definecolor{Green}{rgb}{0,1,0}
\definecolor{Blue}{rgb}{0,0,1}

\begin{document}
\title{Fast and Reliable Primary Frequency Reserves From Refrigerators with Decentralized Stochastic Control}
\author[1]{Evangelos~Vrettos\thanks{vrettos@eeh.ee.ethz.ch}}
\author[2]{Charalampos~Ziras\thanks{chazi@elektro.dtu.dk}}
\author[1]{G\"{o}ran~Andersson\thanks{andersson@eeh.ee.ethz.ch}}
\affil[1]{Power Systems Laboratory, ETH Zurich, Switzerland}
\affil[2]{Department of Electrical Engineering, Technical University of Denmark (DTU), Lyngby, Denmark}
\maketitle

\begin{abstract}
Due to increasing shares of renewable energy sources, more frequency reserves are required to maintain power system stability. In this paper, we present a decentralized control scheme that allows a large aggregation of refrigerators to provide Primary Frequency Control (PFC) reserves to the grid based on local frequency measurements and without communication.

The control is based on stochastic switching of refrigerators depending on the frequency deviation. We develop methods to account for typical lockout constraints of compressors and increased power consumption during the startup phase. In addition, we propose a procedure to dynamically reset the thermostat temperature limits in order to provide reliable PFC reserves, as well as a corrective temperature feedback loop to build robustness to biased frequency deviations. Furthermore, we introduce an additional randomization layer in the controller to account for thermostat resolution limitations, and finally, we modify the control design to account for refrigerator door openings.

Extensive simulations with actual frequency signal data and with different aggregation sizes, load characteristics, and control parameters, demonstrate that the proposed controller outperforms a relevant state-of-the-art controller.
\end{abstract}

\newpage

\begin{acronym}[t]
\acro{CDF}{Cumulative Distribution Function}
\acro{COP}{Coefficient of Performance}
\acro{PDF}{Probability Distribution Function}
\acro{PFC}{Primary Frequency Control}
\acro{SFC}{Secondary Frequency Control}
\acro{RES}{Renewable Energy Sources}
\acro{RMSE}{Root Mean Square Error}
\acro{MAPE}{Mean Absolute Percentage Error}
\acro{MHSE}{Moving Horizon State Estimator}
\acro{TCL}{Thermostatically Controlled Load}
\acro{TSO}{Transmission System Operator}
\end{acronym}

\newpage

\begin{table}[t]
\scriptsize
\renewcommand{\arraystretch}{1.1}
\caption{Nomenclature: symbols}
\centering
\begin{tabular}{|ccc|}
\hline
\textbf{Symbol} & \textbf{Unit} & \textbf{Description}\\
$C$ & kJ/$^\circ$C & Refrigerator thermal capacitance\\
$e_\text{r}$ & $\%$ & Instantaneous reserve error\\
$e_\text{b,mape}$ & $\%$ & Baseline MAPE\\
$e_\text{r,mape}$ & $\%$ & Reserve MAPE\\
$e_{\text{p,rmse}}$ & - & Switching probability RMSE\\
$D$ & - & Duty cycle\\
$D^\text{a}$ & - & Actual duty cycle (with PFC)\\
$D^\text{n}$ & - & Nominal duty cycle (without PFC)\\
$D^\text{d}$ & - & Desired duty cycle (with PFC)\\
$D^\text{r}$ & - & PFC reserve capacity in terms of duty cycle\\
$E_\text{cl}^\text{d}$ & kJ & Daily energy consumption without door openings\\
$E_\text{op}^\text{d}$ & kJ & Daily energy consumption with door openings\\
$F_\text{off}$ & - & CDF of the lock-off time in the aggregation\\
$F_\text{on}$ & - & CDF of the lock-on time in the aggregation\\
$K_{\text{c}}$ & - & Corrective temperature gain\\
$K_{\text{r},t}$ & $^\circ$C & Resetting factor for the thermostat limits\\
$l_\text{off}$ & - & Lock-off event as a random variable\\
$l_\text{on}$ & - & Lock-on event as a random variable\\
$L_\text{off}^\text{st}$ & - & Steady-state fraction of refrigerators locked-off\\
$L_\text{on}^\text{st}$ & - & Steady-state fraction of refrigerators locked-on\\
$L_{\textrm{off},t}^{\textrm{tr}}$ & - & Transient fraction of refrigerators locked-off\\
$L_{\textrm{on},t}^{\textrm{tr}}$ & - & Transient fraction of refrigerators locked-on\\
$m$ & - & Refrigerator on/off state\\
\hline
\end{tabular}
\label{tab:nomenclature1_prim_fridges}
\end{table}

\begin{table}[t]
\scriptsize
\renewcommand{\arraystretch}{1.1}
\caption{Nomenclature: symbols (continued)}
\centering
\begin{tabular}{|ccc|}
\hline
\textbf{Symbol} & \textbf{Unit} & \textbf{Description}\\
$N_\text{ev}$ & s & Duration of a step frequency deviation event\\
$N_\text{on}$ & s & Maximum lock-on time within the aggregation\\
$N_\text{off}$ & s & Maximum lock-off time within the aggregation\\
$N_\text{r}$ & - & Size of refrigerator aggregation\\
$N_\text{rec}$ & s & Recovery period after a step frequency deviation\\
$N_\text{r,a}$ & - & Number of refrigerators activated for PFC\\
$N_\text{s}$ & s & Duration of startup dynamics\\
$N_\text{sim}$ & s & Simulation period\\
$N_\delta$ & s & Half-period of the frequency deviation signal\\
$P_\text{agg}$ & W & Actual aggregate electric power\\
$P_\text{b}$ & W & Uncontrolled aggregate electric power\\
$P_\text{d}$ & W & Desired aggregate electric power (with PFC)\\
$P_\text{res}$ & W & Reserve capacity of the aggregation\\
$P_\text{n}$ & W & Nominal refrigerator electric power\\
$q$ & - & Switching probability of the benchmark controller\\
$R$ & $^\circ$C/kW & Refrigerator's thermal resistance with closed door\\
$R_\text{op}$ & $^\circ$C/kW & Refrigerator's thermal resistance with open door\\
$r$ & - & Normalization factor for the resetting factor\\
$S_\text{off}$ & - & Survival function of the lock-off time in the aggregation\\
$S_\text{on}$ & - & Survival function of the lock-on time in the aggregation\\
$S_\text{u}$ & - & Refrigerator startup dynamics profile\\
$s_{\textrm{net},t}$ & - & Average net switching rate in the aggregation\\
$s_{\textrm{off},t}$ & - & Average rate at which refrigerators switch off\\
$s_{\textrm{on},t}$ & - & Average rate at which refrigerators switch on\\
\hline
\end{tabular}
\label{tab:nomenclature2_prim_fridges}
\end{table}

\begin{table}[t]
\scriptsize
\renewcommand{\arraystretch}{1.1}
\caption{Nomenclature: symbols (continued)}
\centering
\begin{tabular}{|ccc|}
\hline
\textbf{Symbol} & \textbf{Unit} & \textbf{Description}\\
$T$ & $^\circ$C & Refrigerator temperature\\
$T_\text{a}$ & $^\circ$C & Ambient temperature\\
$T_\text{max}$ & $^\circ$C & Higher deadband limit of thermostat\\
$T_\text{min}$ & $^\circ$C & Lower deadband limit of thermostat\\
$\bar{T}_\text{nom}$ & $^\circ$C & Nominal mean temperature (without PFC)\\
$T_\text{set}$ & $^\circ$C & Thermostat setpoint $\left[T_\text{set}=0.5 \cdot (T_\text{max}+T_\text{min})\right]$\\
$\dot{T}_{\textrm{d}}$ & $^\circ$C/s & Average temperature decrease rate\\
$\dot{T}_{\textrm{i}}$ & $^\circ$C/s & Average temperature increase rate\\
$\bar{T}_{\textrm{cl}}$ & $^\circ$C & Average refrigerator temperature with closed door\\
$\bar{T}_{\textrm{op}}$ & $^\circ$C & Average refrigerator temperature with open door\\
$t_\text{off}$ & s & Duration of refrigerator's off cycle\\
$t_\text{off}^l$ & s & Refrigerator lock-off time\\
$t_\text{on}$ & s & Duration of refrigerator's on cycle\\
$t_\text{on}^l$ & s & Refrigerator lock-on time\\
$u$ & - & Peak power factor for startup dynamics modeling\\
$w$ & $^\circ$C/s & Noise term for refrigerator's external disturbances\\
$x$ & - & Fraction of switched loads for PFC under startup dynamics\\
$z$ & - & Thermostat resetting event as a Bernoulli random variable\\
\hline
\end{tabular}
\label{tab:nomenclature3_prim_fridges}
\end{table}

\begin{table}[t]
\scriptsize
\renewcommand{\arraystretch}{1.1}
\caption{Nomenclature: Greek letters}
\centering
\begin{tabular}{|ccc|}
\hline
\textbf{Symbol} & \textbf{Unit} & \textbf{Description}\\
$\alpha$ & $1$/s & Thermal parameter of continuous-time model\\
$\beta$ & $^\circ$C/kJ & Thermal parameter of continuous-time model\\
$\gamma$ & $^\circ$C/Hz & Auxiliary variable for PFC\\
$\Delta{D}$ & - & Duty cycle change\\
$\Delta{D}_t^\text{op}$ & - & Increase in duty cycle due to door openings\\
$\zeta$ & - & Auxiliary variable for thermostat resetting events\\
$\Delta{f}$ & Hz & Frequency deviation\\
$\Delta{T}$ & $^\circ$C & Width of thermostat's deadband\\
$\Delta{T}_\text{b}$ & $^\circ$C & Bound for refrigerator's thermostat limits deviation\\
$\Delta{T}_\text{lim}$ & $^\circ$C & Change in refrigerator's thermostat limits\\
$\Delta{T}_\text{res}$ & $^\circ$C & Minimum resolution of refrigerator's thermostat\\
$\Delta{t}$ & s & Discretization time step\\
$\delta$ & Hz & Magnitude of frequency deviation signal\\
$\varepsilon$ & $^\circ$C & Tolerance for temperature deviation during a frequency event\\
$\epsilon$ & $^\circ$C & Tolerance for temperature deviation after a frequency event\\
$\eta$ & - & Coefficient of performance\\
$\lambda$ & - & Auxiliary variable equal to $1-K_{\text{c}}$\\
$\mu_\text{d}$ & s & Mean duration of a door opening event\\
$\mu_\text{op}$ & - & Average number of door openings per day\\
$\nu$ & - & Parameter for modeling of startup dynamics\\
$\xi$ & - & Normalized energy consumption increase due to door openings\\
$\rho$ & - & Switching probability of the proposed controller\\
$\varrho$ & - & Switching probability for resetting of thermostat limits\\
$\sigma_\text{d}$ & s & Standard deviation of the duration of a door opening event\\
$\sigma_\text{op}$ & - & Standard deviation of the number of door openings per day\\
$\chi$ & - & Switching action of a refrigerator (Bernoulli random variable)\\
\hline
\end{tabular}
\label{tab:nomenclature4_prim_fridges}
\end{table}

\begin{table}[t]
\scriptsize
\renewcommand{\arraystretch}{1.1}
\caption{Nomenclature: probability operators}
\centering
\begin{tabular}{|cc|}
\hline
\textbf{Operator} & \textbf{Description}\\
$\Mean\left[\cdot\right]$ & Expected value of a random variable\\
$\Normal$ & Normal probability distribution\\
$\Prob\left[\cdot\right]$ & Probability of a random variable\\
$\mathrm{SD}\left[\cdot\right]$ & Standard deviation of a random variable\\
$\Unif$ & Uniform probability distribution\\
$\mathrm{Var}\left[\cdot\right]$ & Variance of a random variable\\
\hline
\end{tabular}
\label{tab:nomenclature5_prim_fridges}
\end{table}

\begin{table}[t]
\scriptsize
\renewcommand{\arraystretch}{1.1}
\caption{Nomenclature: subscripts and superscripts}
\centering
\begin{tabular}{|cc|}
\hline
\textbf{Subscript/Superscript} & \textbf{Description}\\
$i$ & Refrigerator index within an aggregation\\
min/max & Minimum/maximum value of a variable or parameter\\
$k, t$ & Time indexes\\
$\bar{(\cdot)}$ & Average value of a variable or parameter\\
$\hat{(\cdot)}$ & Estimated value of a variable or parameter\\
$[\cdot]_+$ & Caps a variable to non-negative values\\
\hline
\end{tabular}
\label{tab:nomenclature6_prim_fridges}
\end{table}

\section{Introduction}
\subsection{Motivation and Related Work}
It is expected that the increasing shares of intermittent \ac{RES} will increase the need for ancillary services in power systems \cite{makarov_operational_2009}. Since \ac{RES} displace conventional power plants, there is a growing interest in exploiting the flexibility of demand-side resources to provide ancillary services, namely frequency and voltage control \cite{callaway_achieving_2011}. \acp{TCL} with thermal inertia such as space and water heaters, air conditioners and refrigerators are well suited for provision of ancillary services, because their consumption can be shifted in time without user discomfort.

The idea of controlling loads to provide ancillary services is not new. In $1980$, \cite{Schweppe1980} proposed a simple frequency-responsive controller to adjust the temperature limits of \acp{TCL} based on frequency measurements to provide a governor-type action. Much of the early work focused on modeling populations of TCLs based on first principles \cite{chong1979statistical,malhame1985electric,Constantopoulos1991}. For example, \cite{malhame1985electric} proposed a well-known model consisting of a set of Fokker-Planck equations that describe the \ac{PDF} of temperature in a TCL aggregation.

The first line of research on load control has focused on centralized control of TCL aggregations to provide load following or frequency regulation. References \cite{callaway2009tapping,Perfumo2012,kundu2011modeling} developed controllers to provide such services via thermostat setpoint manipulation, and showed that linear models are sufficient to describe the aggregate transient dynamics of TCLs under thermostat setpoint control.

Other researchers developed control strategies based on probabilistic switching, where a fraction of the devices switches stochastically to provide the desired service. Reference~\cite{Keep2011switches} showed how bi-directional control can be achieved via broadcasting a signal that only switches off devices. Broadcast controllers were proposed in \cite{Zhang2013Aggregated} and \cite{TotuSwitching2014} to allow the aggregate power of a population of air conditioners and refrigerators to track a power reference signal. Reference~\cite{moya2014hierarchical} proposed a hierarchical load control framework, where the upper layer is centralized and computes optimal control gains for loads at different buses, and the lower layer is decentralized with switching probabilities that depend on the control gains.

The above works as well as other relevant approaches that rely either on probabilistic switching or rule-based controllers (e.g., \cite{koch_modeling_2011,kondoh_evaluation_2011,Vrettos2012}) assume access to TCL state measurements (temperatures and on/off states). More recent works have proposed state estimation methods to reduce the needs for real-time communication. For example, \cite{mathieu_state_2013} used a Kalman Filter to estimate the state of a Markov model for a TCL aggregation, whereas \cite{kara_moving_2013} used a similar modeling approach and a \ac{MHSE}. Reference~\cite{borsche_minimizing_2013} developed a particle filter, whereas \cite{MouraCDC2013} proposed a state estimator based on partial differential equations. The authors of \cite{VrettosPSCC2014} developed an \ac{MHSE} to reconstruct the states of individual TCLs using aggregate power measurements from substations and low-frequency TCL state measurements from smart meters.

The second line of research investigated how TCLs can provide frequency support in a decentralized way based on local frequency measurements. References~\cite{FrerisStabilization2007} and \cite{xu2011demand} developed deterministic approaches that rely on a frequency-dependent temperature deadband to provide \ac{PFC}. Reference~\cite{molina2011decentralized} introduced a more generic approach that is based not only on frequency deviation measurements, but also on their evolution over time. A simple rule-based controller with delays was used in \cite{samarakoon2012investigation} in an experimental demonstration to provide PFC with residential appliances.

Such deterministic approaches provide an effective initial response and reduce the frequency excursions after a sudden disturbance.Despite their simplicity, these approaches have an important limitation: they cause rebound effects and tend to synchronize the on/off cycles of individual devices, which might introduce non-decreasing oscillations in frequency \cite{Angeli2012Fridges,kremers2013emergent}.

In order to de-synchronize the responses of TCLs while providing PFC, a few randomized approaches have been proposed in the literature. For example, a control approach with random on/off frequency thresholds and turn-on delay times was used in an experimental demonstration in \cite{PNNL2007}. Reference~\cite{moghadam2014distributed} proposed a randomized load control scheme where each load monitors the frequency over random time instances and responds according to a simple control policy based on frequency thresholds. 

Reference~\cite{Angeli2012Fridges} proposed an interesting stochastic approach to control a population of refrigerators. Each device is modeled as a Markov-jump linear system with transition probabilities between the on and off states that depend indirectly on the frequency deviation. The authors derived closed form expressions for the mean value and the variance of the average temperature within the population, and proved that the closed-loop system is asymptotically stable. However, the control adjusts the properties of the steady state temperature distribution, and therefore results in slow responses that are unacceptable for PFC.

Another relevant decentralized stochastic controller was developed in \cite{tindemans2015decentralized}. The main innovation of this work is that the average heating rate of the population is used as a control variable. The formulations of \cite{tindemans2015decentralized} result in control laws for individual appliances, namely temperature limit changes and switching rates, and allow an aggregator to estimate the available flexibility from the aggregation. However, the control approach of \cite{tindemans2015decentralized} is able to respond to demand reduction requests but not demand increase requests, and thus it is not appropriate for PFC.

More related to our work is \cite{borsche2015new} where a decentralized stochastic controller based on probabilistic switching is proposed for refrigerators to provide PFC. This controller results in fast responses that are suitable for PFC, but it has a number of limitations that will be discussed in Section~\ref{contribution}.

Some papers, for example \cite{biegel2013primary,weckx2013optimal,zhao2012frequency,zhao2014design}, proposed methods to provide PFC at minimum cost or maximum social welfare. Note that we do not address this topic and, therefore, the proposed controller is complementary to the work presented in \cite{biegel2013primary,weckx2013optimal,zhao2012frequency,zhao2014design}. Reference \cite{biegel2013primary} proposed an optimization problem formulation to allocate trigger frequencies to on/off devices such that they collectively provide PFC reserves at minimum cost. However, \cite{biegel2013primary} used a simplified model for on/off devices, which does not model TCLs in full detail. Reference~\cite{weckx2013optimal} proposed a frequency support method for residential loads based on utility functions to maximize social welfare. The authors of \cite{zhao2012frequency,zhao2014design} formulated a load control optimization problem and developed synchronous and asynchronous algorithms to solve it in a decentralized way. References \cite{zhao2012frequency,zhao2014design} used simplified aggregate load models for each bus of the transmission system, but did not propose methods to coordinate a population of on/off devices.

\subsection{Contribution and Organization of this Paper} \label{contribution}
In this paper, we develop a new decentralized stochastic method to provide accurate and reliable PFC reserves from an aggregation of refrigerators.\footnote{Although we consider residential refrigerators, the proposed approach can be applied to other TCLs with compressors and similar dynamics. The main advantage of residential refrigerators over other TCLs is that they continuously cycle on and off and that their power consumption has very little seasonality (due to regulated room temperature). Therefore, the potential for PFC reserve is relatively constant all year round. Refrigerators are more well-suited for PFC than secondary and tertiary control, due to their limited energy capacity.} In contrast to the centralized approaches of \cite{koch_modeling_2011,Keep2011switches,Vrettos2012,Zhang2013Aggregated,TotuSwitching2014,moya2014hierarchical}, the proposed controller relies on switching probabilities that are computed based on local frequency measurements, and therefore there is no need for real-time communication. Our work is similar to \cite{Angeli2012Fridges,moghadam2014distributed,tindemans2015decentralized} in that it randomizes the reaction of individual loads to avoid synchronization. The main advantage of probabilistic switching is a much faster PFC response compared with \cite{moghadam2014distributed} where the time instance of the reaction of each load is randomized, and \cite{Angeli2012Fridges,tindemans2015decentralized} where the devices' switching rates depend on the desired evolution of mean temperature.

A similar probabilistic switching approach for PFC with refrigerators was proposed in \cite{borsche2015new}. This approach relies on four main assumptions: (i) each refrigerator consumes a constant amount of power at the on state; (ii) each refrigerator can cycle on and off arbitrarily often; (iii) the frequency deviation signal is zero-mean; and (iv) the refrigerators are not exposed to door openings. Assumptions (i), (ii) and (iv) were also made in \cite{tindemans2015decentralized}.

In practice, the power consumption of a refrigerator is typically higher at the beginning of each on cycle (the so-called startup dynamics), its compressor has lockout constraints and cannot turn on and off arbitrarily often, and its temperature is affected by door openings. Furthermore, the frequency deviation might be significantly biased to the positive or negative direction for a few hours. Since these assumptions do not hold in practice, the method of \cite{borsche2015new} might fail to provide accurate PFC reserves in a realistic environment, and thus its applicability is questionable.

Our main contribution is a new controller that works in a realistic environment without the limiting assumptions (i) - (iv).\footnote{This paper extends, improves and formalizes the methods presented in our preliminary work \cite{Ziras2015GM}.} Apart from the basic probabilistic switching formulation which is similar to that of \cite{borsche2015new}, the proposed controller consists of several novel components that allow us to drop assumptions (i) - (iv). Furthermore, we show that combining the probabilistic switching with deterministic control of thermostat temperature limits improves the control performance. Nevertheless, if the frequency deviation signal is biased, the controller's performance will deteriorate because the mean temperature of the population will deviate significantly from the nominal value. For this purpose, we develop an additional control loop to regulate the mean temperature of the population, and derive analytical upper and lower bounds for the control gain. In addition, we introduce another randomization layer in the controller to compensate for the typically limited thermostat resolution. We analytically show that the side-effect is an increased variance of the thermostat limit changes which, nevertheless, can be bounded. In comparison with the controller of \cite{borsche2015new}, the proposed controller achieves a significant improvement in the continuous operation of a refrigerator aggregation under PFC provision, as well as a small improvement in dynamic response.

The remainder of this paper is organized as follows. Section~\ref{model_and_basic_prob_switching} introduces the refrigerator model and probabilistic switching control approach. In Section~\ref{startup_and_lock}, we extend the probabilistic switching to account for startup dynamics and lockout constraints. Section~\ref{resetting_and_corrective} discusses the advantages and implications of controlling the thermostat temperature limits, whereas in Section~\ref{door_openings} we integrate refrigerator door openings in the control scheme. In Section~\ref{parameters_and_metrics}, we introduce the model parameters and performance metrics, which are used in the simulations and sensitivity analysis of Section~\ref{sensitivity}. In Section~\ref{robustness_analysis}, we show the controller's robustness to biased frequency deviations, limited thermostat resolution, and door openings. Implementation issues are discussed in Section~\ref{implementation_issues}, whereas Section~\ref{conclusions} concludes.

\section{Modeling and Probabilistic Switching} \label{model_and_basic_prob_switching}
\subsection{Refrigerator Model}
Consider a refrigerator $i$ without freezer controlled by an on/off hysteresis controller. We use a standard first order differential equation to model the lumped temperature $T_i(t)$ of the refrigerator including the air and solid mass (refrigerator compartment and content) \cite{Constantopoulos1991,Angeli2012Fridges}
\begin{align}
\dot{T}_i(t) &= \alpha_i \big[T_{\textrm{a},i}-T_i(t)\big] + w_i(t) &\text{when OFF} \label{eq:ode_off} \\
\dot{T}_i(t) &= \alpha_i \big[T_{\textrm{a},i}-T_i(t)\big] - \beta_i P_{\textrm{n},i} + w_i(t) &\text{when ON}, \label{eq:ode_on}
\end{align}
where $\alpha = 1/RC$, $\beta = \eta/C$, $C$ is the thermal capacitance, $R$ is the thermal resistance, $T_\textrm{a}$ is the room temperature, $\eta$ is the \ac{COP}, and $P_{\textrm{n}}$ is the nominal power. The noise term $w(t)$ aggregates the effect of external disturbances, e.g., door openings, changes in food content, and variations of $T_{\textrm{a}}$. In Section~\ref{door_openings}, we model the effect of door openings that are the dominant disturbances.

The model is discretized with a discretization time step $\Delta{t}=1~\text{sec}$. Let $T_{\text{min},i}$, $T_{\text{max},i}$ denote the deadband limits of the thermostat's hysteresis controller and $m_{i,t} \in \{0,1\}$ denote the compressor's on/off state. Using the discrete time model we can derive the duration of the on and off cycles ($t_{\textrm{on}}$ and $t_{\textrm{off}}$)
\begin{align}
t_{\textrm{on},i} &= R_i C_i ~ \ln\bigg[\frac{T_{\textrm{max},i}-T_{\textrm{a},i}+\eta_i R_i P_{\textrm{n},i}}{T_{\textrm{min},i}-T_{\textrm{a},i}+\eta_i R_i P_{\textrm{n},i}}\bigg] \label{on_time_est} \\
t_{\textrm{off},i} &= R_iC_i ~ \ln \bigg[ \frac{T_{\textrm{min},i}-T_{\textrm{a},i}}{T_{\textrm{max},i}-T_{\textrm{a},i}} \bigg]~, \label{off_time_est}
\end{align}
and the device's duty cycle $D_i = t_{\textrm{on},i}/(t_{\textrm{on},i}+t_{\textrm{off},i})$.

\subsection{Aggregation Model}
Consider a population of $N_\textrm{r}$ refrigerators providing \ac{PFC}. For control purposes we are interested in the aggregate power of the population $P_{\text{agg},t}$, which depends on the aggregate duty cycle $D_t$, i.e., the fraction of loads that are at the on state. If $N_\textrm{r}$ is sufficiently large, $D_t$ and $P_{\text{agg},t}$ are given by
\begin{align}
D_t = (1/N_\textrm{r}) \cdot \sum\nolimits_i m_{i,t}~,~~P_{\text{agg},t} \approx N_\textrm{r} D_t \bar{P}_{\textrm{n}}~,
\end{align}
where $i \in [1,N_\textrm{r}]$ is the refrigerator index and $\bar{\cdot}$ denotes the mean value of a parameter in the aggregation. Without \ac{PFC} and if $w_i(t)=0$~$\forall i$, the duty cycle will be approximately constant and equal to a nominal value $D^{\textrm{n}}$, whereas the aggregate power will be approximately equal to $P_{\text{agg},t} \approx N_\textrm{r} D^{\textrm{n}} \bar{P}_{\textrm{n}}$.

\subsection{Primary Frequency Control with Probabilistic Switching} \label{prob_switching}
\ac{PFC} reserves can be provided by perturbing the population's duty cycle around $D^{\textrm{n}}$ proportionally to the frequency deviation $\Delta{f}_t$ with a gain $D^{\textrm{r}}$, which corresponds to the \ac{PFC} reserve capacity (the reserve capacity increases with the value of $D^{\textrm{r}}$). Therefore, in order to respond to a frequency deviation $\Delta{f}_t$ the duty cycle should be controlled to the desired value
\begin{align}
D^{\textrm{d}}_t = D^{\textrm{n}} + D^{\textrm{r}} \cdot \frac{\Delta{f}_t}{\Delta{f}_{\textrm{max}}}~, \label{desired_duty_cycle}
\end{align}
where $\Delta{f}_{\textrm{max}}$ is the frequency deviation for full \ac{PFC} reserve activation, e.g., $0.2$~Hz in continental Europe.

Let us denote by $\Delta{D}_t = D^{\textrm{d}}_t - D^{\textrm{d}}_{t-1}$ the change in desired duty cycle between two consecutive time steps. Note that $\Delta{D}_t$ depends on $\Delta{f}_t$ and $\Delta{f}_{t-1}$, which can be measured locally by each refrigerator, and so the duty cycle change can be achieved in a decentralized way with probabilistic switching. Each device calculates the switching probability
\begin{equation}
q_t = \begin{cases}
        \frac{\Delta{D}_t}{1-D^{\textrm{d}}_{t-1}}, & \textrm{if} \hspace{2mm} \Delta{D}_t \geq 0~~\text{(switch on)}\\					 \label{eq:k}
        -\frac{\Delta{D}_t}{D^{\textrm{d}}_{t-1}}, & \textrm{if} \hspace{2mm} \Delta{D}_t<0~~\text{(switch off)}~,
        \end{cases}
\end{equation}
and generates a random number uniformly distributed in $[0,1]$. If the number is smaller than $q_t$, the device will switch to contribute to PFC. The first line of \eqref{eq:k} is the probability for devices that are off to switch on; thus, the desired duty cycle change $\Delta{D}_t$ is normalized by the portion of devices that are currently off ($1-D^{\textrm{d}}_{t-1}$). The second line of \eqref{eq:k} is the probability for devices that are on to switch off, and $\Delta{D}_t$ is normalized by the portion of devices that are currently on ($D^{\textrm{d}}_{t-1}$).

If each refrigerator responds to the switching probability \eqref{eq:k}, the fraction of switched devices will be very close to $\Delta{D}_t$ due to the law of large numbers. Therefore, the population will collectively provide the required reserve without any real-time communication between the loads.

\subsection{Offline Communication and Coordination}\label{offline_coordination}
Although no real-time communication is required to provide PFC, the refrigerators need some offline coordination. The main reason for this is to fix the value of $D^{\textrm{r}}$ that reflects the PFC reserve capacity, which is typically determined in a reserve market on a weekly or daily basis.\footnote{{A refrigerator aggregation can provide PFC in parallel to other reserve providers (e.g., generators or demand-side resources) without any coordination with them other than the market-based determination of reserve capacities.}} Furthermore, as discussed in the following sections, the control design relies on mean values of several parameters in the population (e.g., $\bar{P}_{\textrm{n}}$). Therefore, if the statistical properties of the population change over time, the respective controller parameters should be updated. A communication channel with low data transfer can be used to (infrequently) coordinate the refrigerators by communicating the updated controller parameters and $D^{\textrm{r}}$ values.

\section{Startup Dynamics and Lockout Constraints} \label{startup_and_lock}
\subsection{Refrigerator Startup Dynamics} \label{startup_dyn}
Typically, a refrigerator's compressor consumes more power during the startup phase due to a smaller \ac{COP}, higher evaporator temperature, and higher motor current \cite{Goncalves2011Fridges}. Since this additional power is up to $25\%$ of the nominal power $P_{\textrm{n}}$, it is important to consider it in \ac{PFC} to avoid overshoots and large reserve errors.

We model the startup dynamics of refrigerator $i$ using a peak power factor $u_i > 0$ and the duration of the startup phase $N_{\textrm{s},i}$. We assume that the power becomes immediately equal to \mbox{$(1+u_i) P_{\textrm{n},i}$} when the refrigerator switches on, and then it linearly decreases until it becomes equal to $P_{\textrm{n},i}$ after $N_{\textrm{s},i}$ seconds. Thus, the refrigerator power evolves according to
\begin{align}
P_{i,t} = P_{\textrm{n},i} \hspace{-0.05cm} \left[1+S_{\text{u},i}(t)\right] = P_{\textrm{n},i} \hspace{-0.05cm} \Bigg[1+ u_i \bigg[1-\frac{t}{N_{\textrm{s},i}}\bigg]_+\Bigg]~, \label{peak_dyn}
\end{align}
where $S_{\text{u},i}$ is the startup power profile and $[\cdot]_+$ caps its argument to positive values, i.e., $S_{\text{u},i}(t)=0$ for $t \geq N_{\textrm{s},i}$.

If the power consumption of each refrigerator was constant over time, $\Delta{D}_t$ in \eqref{eq:k} would be the fraction of loads that need to switch, and it would depend only on $\Delta{f}_t$ and $\Delta{f}_{t-1}$. Since the startup dynamics affect the future aggregate power, the fraction of loads that need to switch at the current time step depends additionally on the number of loads that switched during the previous reserve activations. Therefore, the fraction of loads to be switched -- denoted by $x_t$ -- is in general different from $\Delta{D}_t$.

\begin{figure}[t]
\centering
\includegraphics[width=\textwidth]{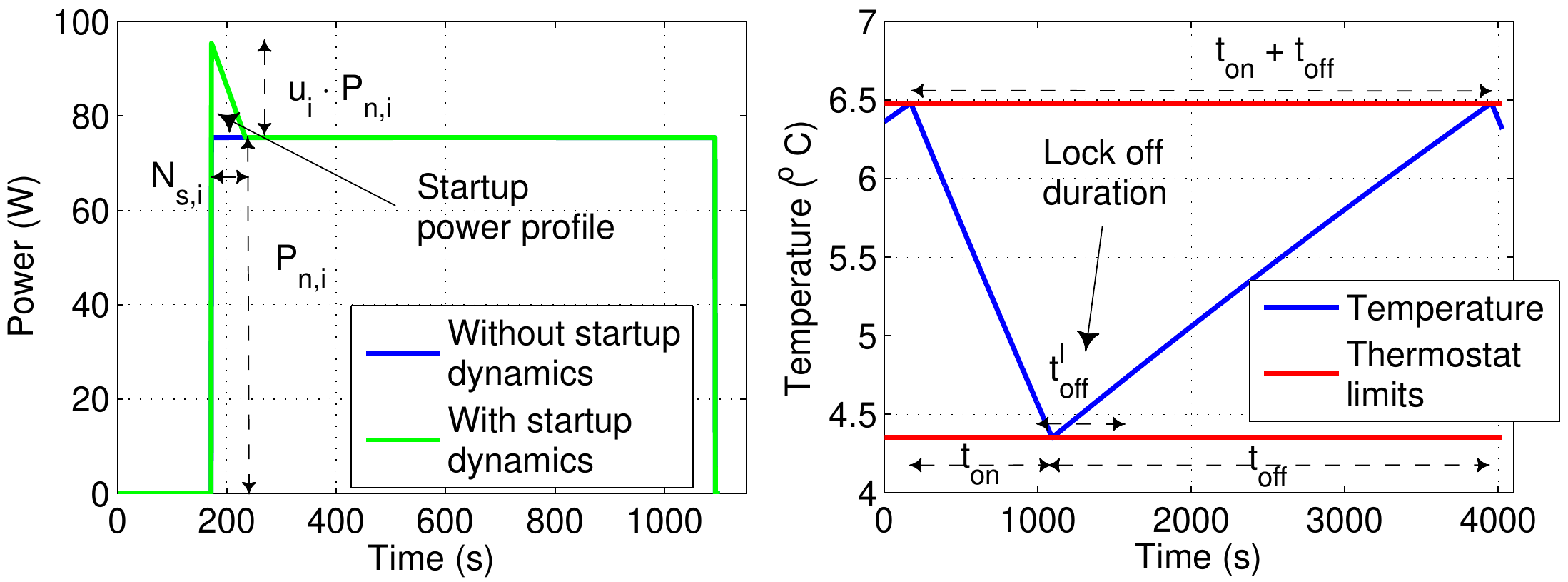}
\caption{Left: Effect of startup dynamics on refrigerator power. Right: Evolution of refrigerator temperature during an on/off cycle.} \label{fig:startup_and_lockout}
\end{figure}

For control design we assume to know only the average values $\bar{P}_\text{n}=\Mean\left[P_{\text{n},i}\right]$, $\bar{u}=\Mean\left[u_i\right]$ and $\bar{N}_{\textrm{s}}=\Mean\left[N_{\textrm{s},i}\right]$ in the population\footnote{The average values can be computed from available statistics on domestic refrigerators. Typical values are given in Table~\ref{table:general_fridge_params}.}, and define the average startup duty cycle profile
\begin{align}
\bar{S}_{\text{u}}(t) = \bar{u} \cdot \left[1-\frac{t}{\bar{N}_{\textrm{s}}}\right]_+ ~, \label{Su_peak_dyn}
\end{align}
which we use to recursively calculate $x_t$ with
\begin{align}
x_t &= \frac{1}{1+\nu_t} \hspace{-0.04cm} \bigg[D^{\textrm{d}}_t-D^{\textrm{a}}_{t-1}- \hspace{-0.09cm} \sum_{i=t-N_{\textrm{s}}+1}^{t-1}x_i \bar{S}_{\textrm{u}}(t-i)\bigg] \label{recurX} \\
D^{\textrm{a}}_t &= D^{\textrm{a}}_{t-1} + x_t~. \label{XandD}
\end{align}
The variable $D^{\textrm{a}}_t$ is the actual duty cycle of the population, whereas $\bar{S}_{\text{u}}(t)$ models the startup dynamics as the triangle shown in the left plot of Fig.~\ref{fig:startup_and_lockout}. For $t=0$, $\bar{S}_{\textrm{u}}(0)$ corresponds to the normalized overshoot in power consumption $\bar{u}$. For example, if the maximum power consumption of refrigerators at the startup phase is $25\%$ higher than the nominal power, then $\bar{S}_{\textrm{u}}(0)=0.25$. For $0 < t \leq N_\text{s}$, the variable $\bar{S}_{\textrm{u}}(t)$ models the decay of the originally higher power to the nominal power.

The term $1/(1+\nu_t)$ in \eqref{recurX} reduces the fraction of switched loads $x_t$ as much as the expected overshoot due to the current startup dynamics. Note that $\nu_t=\bar{S}_{\textrm{u}}(0)$ if the right hand side of \eqref{recurX} is non-negative, whereas $\nu_t=0$ if it is negative. This differentiation is necessary since there exist no dynamics at shutdown. On the other hand, the summation term of \eqref{recurX} recursively accounts for the anticipated power overshoot due to the past startup dynamics associated with recent switching actions up to time step $t-1$.

Using \eqref{recurX} and \eqref{XandD} in the controller, we estimate the effect of refrigerators' startup dynamics on the aggregate power and the fraction of switched loads. Since \eqref{recurX} and \eqref{XandD} rely only on the mean values $\bar{P}_\text{n}$, $\bar{u}$ and $\bar{N}_{\textrm{s}}$, the control will not be perfect. However, it is possible to characterize its performance analytically as shown by Proposition~\ref{startup_proposition}.

\begin{proposition} \label{startup_proposition}
Let us denote by $N_\text{s,min}$ and $N_\text{s,max}$ the minimum and maximum startup duration in the population of refrigerators. The estimated aggregate power due to refrigerator startup dynamics obtained using the average values $\bar{P}_\text{n}$, $\bar{u}$, $\bar{N}_\text{s}$ and \eqref{Su_peak_dyn} is an upper bound of the actual aggregate power up to $N_\text{s,min}$ seconds after the reserve activation, irrespective of the probability distribution of $N_{\textrm{s},i}$.

Furthermore, if $N_{\textrm{s},i}$ follows the uniform distribution $N_{\textrm{s},i}$$\scriptstyle\sim$$[N_\text{s,min}, N_\text{s,max}]$, then the estimated aggregate power is an upper bound for up to $t_\text{lim}$ seconds after the reserve activation, whereas it is a lower bound afterwards, where $t_\text{lim}$ is given by
\begin{align}
t_\text{lim} = \frac{N_\text{s,max} \left(N_\text{s,min}+N_\text{s,max}\right)}{3N_\text{s,max}-N_\text{s,min}}~.
\end{align}
\end{proposition}

\begin{proof}
The proof is given in the appendix.
\end{proof}

Another effect of startup dynamics is an increase in the average power of the refrigerator aggregation, both with and without \ac{PFC}. It is straightforward to express the increased aggregate power as a function of the average values of the parameters of startup dynamics ($\bar{P}_\text{n}$, $\bar{u}$, and $\bar{N}_{\textrm{s}}$), which is then used as a baseline for \ac{PFC} evaluation purposes.

\subsection{Compressor Lockout Constraints} \label{lockout_cons}
Lockout times are usually employed to avoid compressor's frequent switching that decreases efficiency and possibly the lifetime \cite{Zhang2013Aggregated}. Let us denote by $t_{\textrm{on}}^{l}$ the lock-on time, i.e., the minimum duration the device must operate after switching on, and by $t_{\textrm{off}}^{l}$ the lock-off time, i.e., the minimum duration the device must remain off after switching off. When a compressor is locked it cannot react to frequency deviations, and so it is important to dynamically estimate the fraction of locked devices to avoid performance degradation or even loss of controllability. In the following, we present a decentralized method to do so.\footnote{The effect of lockout constraints on the refrigerators' aggregate response resembles a dead-band effect. However, a crucial difference is that the lockout effect is not memoryless, but it evolves in time depending on the statistics of lockout times. For this reason, the standard methods in control theory to handle the deadband effect cannot be used to model the effect of lockout constraints \cite{khalil1996nonlinear}.}

\subsubsection{Steady-state Lockout due to Thermostatic Control Actions} \label{steady_state_lock}
Even without \ac{PFC}, i.e., at steady state, at any time step $t$ a fraction $L_{\textrm{on}}^{\textrm{st}}$ of refrigerators is locked at the on state and a fraction $L_{\textrm{off}}^{\textrm{st}}$ is locked at the off state due to thermostatic control actions. Considering the cycle of an individual refrigerator in Fig.~\ref{fig:startup_and_lockout}, the lock event at the off state $l_\text{off}$ can be modeled as the Bernoulli probability distribution
\begin{align} \label{Bernoulli_off}
f^l_\text{off}(l_\text{off}) = \begin{cases}
                        \frac{t_{\textrm{off}}^{l}}{t_{\textrm{on}}+t_{\textrm{off}}}, & \textrm{if} \hspace{2mm} l_\text{off}=1\\			
                        1-\frac{t_{\textrm{off}}^{l}}{t_{\textrm{on}}+t_{\textrm{off}}}, & \textrm{if} \hspace{2mm} l_\text{off}=0~.
                    \end{cases}
\end{align}

The fraction of refrigerators that is locked at the off state is defined and calculated using the properties of the Bernoulli \ac{PDF}
\begin{align}
L_{\textrm{off}}^{\textrm{st}} := \frac{1}{N_\text{r}} \cdot \sum\nolimits_{i=1}^{N_\text{r}} l_{\text{off},i} = \Mean\left[l_\text{off}\right] = \frac{t_{\textrm{off}}^{l}}{t_{\textrm{on}}+t_{\textrm{off}}}~. \label{fraction_lock_off1}
\end{align}
Since $t_{\textrm{off}}^{l}$, $t_{\textrm{on}}$ and $t_{\textrm{off}}$ are random variables, $L_{\textrm{off}}^{\textrm{st}}$ is also a random variable with a so-called ratio distribution. Without knowledge of the \ac{PDF} of $t_{\textrm{off}}^{l}$, $t_{\textrm{on}}$, $t_{\textrm{off}}$, the expectation of $L_{\textrm{off}}^{\textrm{st}}$ can be approximated with a first order Taylor expansion \cite{van2000mean}
\begin{align}
\Mean\left[L_{\textrm{off}}^{\textrm{st}}\right] = \Mean\bigg[\frac{t_{\textrm{off}}^{l}}{t_{\textrm{on}}+t_{\textrm{off}}}\bigg] \approx \frac{\Mean\left[t_{\textrm{off}}^{l}\right]}{\Mean\left[t_{\textrm{on}}\right]+\Mean\left[t_{\textrm{off}}\right]}~. \label{fraction_lock_off2}
\end{align}
Using a similar procedure, one can get the approximation for steady-state lockout at the on state
\begin{align}
\Mean\left[L_{\textrm{on}}^{\textrm{st}}\right] \approx \frac{\Mean\left[t_{\textrm{on}}^{l}\right]}{\Mean\left[t_{\textrm{on}}\right]+\Mean\left[t_{\textrm{off}}\right]}~. \label{fraction_lock_on2}
\end{align}

\subsubsection{Transient Lockout due to Primary Frequency Control}
When providing \ac{PFC} the refrigerators are exposed to additional switching actions that increase the fraction of locked devices. Let us denote by $L_{\textrm{on},t}^{\textrm{tr}}$, $L_{\textrm{off},t}^{\textrm{tr}}$ the time-varying fractions of devices that lock at the on and off states due to \ac{PFC}. Since $L_{\textrm{on},t}^{\textrm{tr}}$ and $L_{\textrm{off},t}^{\textrm{tr}}$ depend on the past reserve activations, they can be estimated recursively given sufficient statistical information. We assume that the \acp{CDF} of the lockout times are known, and we denote them by $F_{\textrm{on}}(t_{\textrm{on}}^{l})$  for lock-on time and $F_{\textrm{off}}(t_{\textrm{off}}^{l})$ for lock-off time. We further define $N_{\textrm{on}} = \max \left[t_{\textrm{on}}^{l}\right]$ and $N_{\textrm{off}} = \max \left[t_{\textrm{off}}^{l}\right]$. Then, $L_{\textrm{on},t}^{\textrm{tr}}$ and $L_{\textrm{off},t}^{\textrm{tr}}$ can be computed as
\begin{align}
L_{\textrm{on},t}^{\textrm{tr}} \hspace{-0.05cm}&=\hspace{-0.05cm} \sum\nolimits_{k=0}^{t-1} c_k x_k S_{\textrm{on}}(t-k),~ c_k = \begin{cases} 1, \hspace{-0.15cm}& \textrm{if} \hspace{1.2mm} x_k\geq0 \vspace{-0.15cm}\\ 0, \hspace{-0.15cm}& \textrm{if} \hspace{1.2mm} x_k<0 \end{cases}  \label{l_ON} \\
L_{\textrm{off},t}^{\textrm{tr}} \hspace{-0.05cm}&=\hspace{-0.05cm} \sum\nolimits_{k=0}^{t-1} d_k x_k S_{\textrm{off}}(t-k),~ d_k = \begin{cases} 0, \hspace{-0.15cm}& \textrm{if} \hspace{1.2mm} x_k\geq0 \vspace{-0.15cm}\\ 1, \hspace{-0.15cm}& \textrm{if} \hspace{1.2mm} x_k<0 \end{cases}, \label{l_OFF}
\end{align}
where $S_{\textrm{on}}=1-F_{\textrm{on}}$ and $S_{\textrm{off}}=1-F_{\textrm{off}}$ are the survival functions of the \acp{CDF}, $S_{\textrm{on}}(t-k)=0$ for $t-k>N_{\textrm{on}}$, and $S_{\textrm{off}}(t-k)=0$ for $t-k>N_{\textrm{off}}$.

The total fraction of locked devices is simply obtained by adding the steady-state and transient contributions
\begin{align}
L_{\textrm{on},t} = \Mean\left[L_{\textrm{on}}^{\textrm{st}}\right]+L_{\textrm{on},t}^{\textrm{tr}},~
L_{\textrm{off},t} = \Mean\left[L_{\textrm{off}}^{\textrm{st}}\right]+L_{\textrm{off},t}^{\textrm{tr}}~. \label{eq:countONOFF}
\end{align}

\subsection{Improved Probabilistic Switching} \label{prob_switching_final}
In Sections~\ref{startup_dyn} and \ref{lockout_cons} we presented methods to account for the startup dynamics and lockout constraints when providing \ac{PFC} with an aggregation of refrigerators. To integrate these methods in the probabilistic switching approach of Section~\ref{prob_switching}, it suffices to apply the switching probability
\begin{equation}
\rho_t = \begin{cases}
        \frac{x_t}{1-D^{\textrm{a}}_{t-1}-L_{\text{off},t-1}}, & \textrm{if} \hspace{2mm} x_t\geq 0~~\text{(switch on)}\\					 \label{eq:rho}
        -\frac{x_t}{D^{\textrm{a}}_{t-1}-L_{\text{on},t-1}}, & \textrm{if} \hspace{2mm} x_t<0~~\text{(switch off)}~,
        \end{cases}
\end{equation}
instead of \eqref{eq:k}. Observe that the total fraction of loads that need to switch ($x_t$) is normalized by the fraction of loads that are available to respond. If $x_t>0$, the loads that can respond are those that are at the off state and unlocked, i.e., the normalization factor is $1-D^{\textrm{a}}_{t-1}-L_{\text{off},t-1}$. On the other hand, if $x_t<0$, the available loads are those that are at the on state and unlocked, i.e., the normalization factor is $D^{\textrm{a}}_{t-1}-L_{\text{on},t-1}$.

In order to further reduce the compressor cycling, the switching probability can be defined as a function of the time elapsed since the latest switching action of each device. However, we do not use a time-dependent probability calculation in our controller, because it tends to synchronize the load cycles and introduce long-term oscillations in aggregate power.

\section{Control of Thermostat Temperature Limits} \label{resetting_and_corrective}
\subsection{Resetting of Thermostat Limits} \label{resetting_thermosta_limits}
The probabilistic switching of Section~\ref{prob_switching_final} provides accurate \ac{PFC} reserves immediately after a frequency deviation. However, if the frequency deviation remains, the accuracy of reserve provision will deteriorate because the refrigerators' aggregate power will decay towards the pre-disturbance value. This is straightforward to verify if the startup dynamics and lockout constraints are neglected.

Let $s_{\textrm{off},t}$ and $s_{\textrm{on},t}$ denote the rates at which refrigerators switch off and on at time step $t$, which we approximate with
\begin{align}
s_{\textrm{off},t}\approx D^{\textrm{d}}_t \cdot \dot{T}_{\textrm{d}}(\bar{T}),~
s_{\textrm{on},t}\approx (1-D^{\textrm{d}}_t) \cdot \dot{T}_{\textrm{i}}(\bar{T})~, \label{son_soff}
\end{align}
where $\dot{T}_{\textrm{d}}(\bar{T})$ and $\dot{T}_{\textrm{i}}(\bar{T})$ are the temperature decrease and increase rates evaluated at the mean temperature, respectively. $\dot{T}_{\textrm{d}}(\bar{T})$ can be obtained from \eqref{eq:ode_on} and $\dot{T}_{\textrm{i}}(\bar{T})$ from \eqref{eq:ode_off}, and for the average values of the parameters used in this paper we get $\dot{T}_{\textrm{d}}(\bar{T})=-0.0026^\circ\textrm{C/s}$ and $\dot{T}_{\textrm{i}}(\bar{T})=0.0009^\circ\textrm{C/s}$. Combining \eqref{son_soff} and \eqref{desired_duty_cycle}, the net switching rate can be expressed as
\begin{align}
&s_{\textrm{net},t} = s_{\textrm{off},t} + s_{\textrm{on},t} \nonumber \\
&\approx D^{\textrm{n}} \dot{T}_{\textrm{d}}(\bar{T}) + (1-D^{\textrm{n}}) \dot{T}_{\textrm{i}}(\bar{T}) + D^{\textrm{r}} \frac{\Delta{f}_t}{\Delta{f}_{\textrm{max}}} \big[\dot{T}_{\textrm{d}}(\bar{T})-\dot{T}_{\textrm{i}}(\bar{T})\big] \nonumber \\ &= D^{\textrm{r}} \cdot (\Delta{f}_t/\Delta{f}_{\textrm{max}}) \cdot \big[\dot{T}_{\textrm{d}}(\bar{T})-\dot{T}_{\textrm{i}}(\bar{T})\big]~. \label{netSwitching}
\end{align}
In \eqref{netSwitching} $s_{\textrm{off}}=D^{\textrm{n}} \cdot \dot{T}_{\textrm{d}}(\bar{T})$ and $s_{\textrm{on}}=(1-D^{\textrm{n}}) \cdot \dot{T}_{\textrm{i}}(\bar{T})$ are the switching rates at steady state and $s_{\textrm{off}} + s_{\textrm{on}} =0$ holds because the aggregate baseline power is approximately constant.

If $\Delta{f}_t>0$, some refrigerators will switch on to provide \ac{PFC} and $s_{\textrm{net},t}$ will become negative because $\dot{T}_{\textrm{d}}(\bar{T})-\dot{T}_{\textrm{i}}(\bar{T})<0$. This means that some refrigerators will start switching off and therefore the aggregate power will start decreasing towards the baseline. Similarly, if $\Delta{f}_t<0$ the initially lower aggregate power will start increasing towards the baseline. Thus, the reserve accuracy will deteriorate over time.

A solution to this problem is to dynamically modify the devices' thermostat limits in order to keep $s_{\textrm{net},t}$ close to zero. For this purpose, we apply the temperature resetting factor
\begin{align} \label{resetting_factor}
K_{\text{r},t}~[^\circ\textrm{C}] = s_{\textrm{net},t} \cdot \Delta{t} = -D^{\textrm{r}} \Delta{t} \bar{\beta} \bar{P}_\text{n} \cdot \frac{\Delta{f}_t}{\Delta{f}_{\textrm{max}}}~,
\end{align}
where $\dot{T}_{\textrm{d}}(\bar{T})-\dot{T}_{\textrm{i}}(\bar{T})=-\bar{\beta} \bar{P}_\text{n}$ from \eqref{eq:ode_off} and \eqref{eq:ode_on}. If $\Delta{f}_t>0$ the resetting factor is negative, which means that both thermostat limits will decrease with a rate equal to $K_{\text{r},t}$ such that the population reaches an equilibrium at an aggregate power higher than the baseline. Similarly, if $\Delta{f}_t<0$ both thermostat limits increase at a rate equal to $K_{\text{r},t}$.

To account for startup dynamics and lockout constraints a few modifications are needed in \eqref{resetting_factor}. The first modification is to substitute the term $D^{\textrm{r}} \cdot (\Delta{f}_t/\Delta{f}_{\textrm{max}})$ with $x_t$ to account for startup dynamics -- see \eqref{XandD}. The second modification is to incorporate the sequential unlocking of refrigerators after a frequency deviation event based on the \acp{CDF} of lockout times $F_{\textrm{on}}(t_{\textrm{on}}^{l})$ and $F_{\textrm{off}}(t_{\textrm{off}}^{l})$.

Consider a positive $\Delta{f}_t$ starting at time step $t=0$ that will induce a positive $x_t$. Instantaneously, some refrigerators will switch on to provide \ac{PFC} and will lock at the on state. Therefore, at $t=0$ fewer devices will be at the off state, $s_{\textrm{on},t}$ will start decreasing, and an instantaneous resetting factor of $-x_t \dot{T}_{\textrm{i}}(\bar{T})$ is needed to keep $s_{\textrm{net},t}$ close to zero. However, as time elapses, some of the refrigerators that switched on will unlock based on $F_{\textrm{on}}(t_{\textrm{on}}^{l})$, and larger changes in the temperature deadband limits are required to maintain the same aggregate power. Therefore, for $\Delta{f}_t>0$ the resetting factor is determined by $x_t [\dot{T}_{\textrm{d}}(\bar{T}) F_{\textrm{on}}(t) - \dot{T}_{\textrm{i}}(\bar{T})] \Delta{t}$ and for $\Delta{f}_t<0$ by $x_t [\dot{T}_{\textrm{d}}(\bar{T}) - \dot{T}_{\textrm{i}}(\bar{T}) F_{\textrm{off}}(t)] \Delta{t}$.

To generalize for multiple events, the resetting factor at time step $t$ should accumulate all contributions from the previous time steps. Therefore, we compute the resetting factor with
\begin{align} \label{resetting_factor2}
K_{\text{r},t} &= r_t \cdot \Delta{t} \cdot \sum\nolimits_{k=0}^{t-1} x_k \cdot \big[c_k \cdot [\dot{T}_{\textrm{d}}(\bar{T}) \cdot F_{\textrm{on}}(t-k) - \dot{T}_{\textrm{i}}(\bar{T})] \nonumber \\ &+ d_k \cdot [\dot{T}_{\textrm{d}}(\bar{T}) - \dot{T}_{\textrm{i}}(\bar{T}) \cdot F_{\textrm{off}}(t-k)]\big]~,
\end{align}
where $c_k$ and $d_k$ are defined in \eqref{l_ON} and \eqref{l_OFF}. Recall that $F_{\textrm{on}}$, $F_{\textrm{off}}$ are \acp{CDF}, and so $F_{\textrm{on}}(t-k)=1$ for $t-k>N_\text{on}$ and $F_{\textrm{off}}(t-k)=1$ for $t-k>N_\text{off}$ in \eqref{resetting_factor2}. Note that the temperature resetting is applied only to the unlocked refrigerators, because only those affect the net switching rate. This is accounted for in \eqref{resetting_factor2} with the normalization factor
\begin{align}
r_t=\frac{1-\Mean\left[L_{\textrm{on}}^{\textrm{st}}\right]-\Mean\left[L_{\textrm{off}}^{\textrm{st}}\right]}{1-L_{\textrm{on},t}-L_{\textrm{off},t}}~.		 \label{eq:r_t} 
\end{align}
Observe that $r_t>1$ as long as $L_{\textrm{on},t}^{\textrm{tr}}\neq 0$ or $L_{\textrm{off},t}^{\textrm{tr}}\neq 0$, whereas $r_t=1$ if $L_{\textrm{on},t}^{\textrm{tr}} = L_{\textrm{off},t}^{\textrm{tr}} = 0$, i.e., at steady state.

By modifying the thermostat limits using \eqref{resetting_factor2} we effectively transform the electric energy surplus or deficit due to the PFC reserve into thermal energy. In this way, we avoid the decay of aggregate power to the pre-disturbance value and achieve constant reserve provision.

\subsection{Corrective Temperature Control} \label{corrective_control_biased_df}
If the frequency deviation $\Delta{f}_t$ is approximately zero-mean, the population's mean temperature does not change significantly and the duty cycle $D_t$ remains close to the nominal value $D^{\text{n}}$. However, in case of biased frequency deviations the resetting factor $K_{\text{r},t}$ results in prolonged reduction or increase of the thermostat limits, which is undesirable for two main reasons. First, it negatively affects the user utility due to very high or very low temperatures. Second, the population's mean temperature changes substantially, which results in steady-state baseline deviations that affect the generation-demand balance on the system level. A simple solution to this problem would be to filter the frequency deviation signal in order to eliminate the bias \cite{borsche2015new}. However, the disadvantage of this approach is that the bias needs to be absorbed by other available \ac{PFC} resources or transferred to secondary frequency control.

We propose a different approach that allows us to handle biased frequency deviations without relying on additional resources. An opposite control action is imposed on the thermostat limits such that they return to the nominal values in the long run. This \emph{``corrective temperature control''} is designed as a proportional feedback controller on the average temperature across the population using a gain $K_\text{c}$. Thus, the temperature limits of refrigerator $i$ evolve in time according to
\begin{align}
T_{\text{min},i,t} &= T_{\text{min},i,t-1} + K_{\text{r},t} - K_\text{c} \cdot (\bar{T}_{t-1}-\bar{T}_{\text{nom}}) \label{Tmin_corrective}\\
T_{\text{max},i,t} &= T_{\text{max},i,t-1} + K_{\text{r},t} - K_\text{c} \cdot (\bar{T}_{t-1}-\bar{T}_{\text{nom}})~, \label{Tmax_corrective}
\end{align}
where $\bar{T}_{t-1}$ is an estimate of the population's average temperature, and $\bar{T}_{\text{nom}}$ is the nominal mean temperature without \ac{PFC}.

\subsubsection{Estimation of Mean Temperature and Duty Cycle} \label{temp_estimation_section}
The mean temperature $\bar{T}_t$ depends on the past temperature resetting and corrective control actions according to
\begin{align} \label{mean_temp_est}
\bar{T}_t \hspace{-0.06cm} = \hspace{-0.06cm} \bar{T}_{t-1} \hspace{-0.06cm} + \hspace{-0.06cm} [ \hspace{-0.0cm} K_{\text{r},t} \hspace{-0.06cm} - \hspace{-0.06cm} K_\text{c} \hspace{-0.06cm} (\bar{T}_{t-1} \hspace{-0.06cm} - \hspace{-0.06cm} \bar{T}_{\text{nom}}) \hspace{-0.0cm} ] \hspace{-0.08cm} \cdot \hspace{-0.08cm} (1 \hspace{-0.08cm} - \hspace{-0.08cm} L _{\textrm{on},t} \hspace{-0.08cm} - \hspace{-0.08cm} L_{\textrm{off},t}).
\end{align}
The term $\left(1-L_{\textrm{on},t}-L_{\textrm{off},t} \right)$ in \eqref{mean_temp_est} is needed because the temperature resetting and corrective control are applied only to the unlocked devices. With reference to \eqref{fraction_lock_off2} and \eqref{fraction_lock_on2}, the fractions of loads that are locked at steady-state depend on the duration of the on and off cycles. However, as the mean temperature changes, the expected values of the on and off cycles also change and are approximated with
\begin{align}
\Mean\left[t_{\text{on},t}\right] &\approx \bar{R} \bar{C} \cdot \ln\bigg[\frac{\bar{T}_{\textrm{max},t-1}-\bar{T}_\textrm{a}+\bar{\eta} \bar{R} \bar{P}_{\textrm{n}}}{\bar{T}_{\textrm{min},t-1}-\bar{T}_\textrm{a}+\bar{\eta} \bar{R} \bar{P}_{\textrm{n}}}\bigg] \label{on_time_est} \\
\Mean\left[t_{\text{off},t}\right] &\approx \bar{R} \bar{C} \cdot \ln \bigg[ \frac{\bar{T}_\textrm{a}-\bar{T}_{\textrm{min},t-1}}{\bar{T}_\textrm{a}-\bar{T}_{\textrm{max},t-1}} \bigg]~, \label{off_time_est}
\end{align}
where $\bar{T}_{\textrm{min},t-1}=\bar{T}_{t-1}-0.5 \cdot \overline{\Delta{T}}$, $\bar{T}_{\textrm{max},t-1}=\bar{T}_{t-1}+0.5 \cdot \overline{\Delta{T}}$, and $\Delta{T}$ is the deadband width.

Using $\bar{T}_t$ it is also possible to keep track of the population's baseline duty cycle (without PFC) with the approximation
\begin{align} \label{duty_cycle_estimate1}
\bar{D}^\text{n}_t \approx \frac{\Mean\left[t_{\text{on},t}\right]}{\Mean\left[t_{\text{on},t}\right]+\Mean\left[t_{\text{off},t}\right]}~.
\end{align}
The actual duty cycle while providing PFC is given by
\begin{align}
D^{\textrm{a}}_t &= D^{\textrm{a}}_{t-1} + x_t + (\bar{D}^\text{n}_t-\bar{D}^\text{n}_{t-1})~, \label{XandD2}
\end{align}
where $x_t$ is the fraction of loads activated for \ac{PFC}, and \mbox{$\bar{D}^\text{n}_t-\bar{D}^\text{n}_{t-1}$} is the change in the baseline duty cycle due to deviation of the mean temperature from its nominal value.

Although the average temperature can be estimated in a decentralized way using \eqref{mean_temp_est}, the estimates could be reset to the actual values periodically based on measurements of refrigerator temperatures transmitted via a low data transfer communication link (see Section~\ref{offline_coordination}). Apart from providing an initial temperature estimate when the aggregation starts providing PFC reserves, a communication link will also enhance the controller's robustness to estimation error.

\subsubsection{Analytical Tuning of the Corrective Temperature Gain} \label{analytical_Kc}
From a user point of view, high $K_\text{c}$ values are preferable to keep the temperature deviations low. However, from a power system point of view, there is a tradeoff to consider when tuning $K_\text{c}$. If the gain is very small, prolonged temperature deviations will change the aggregation's baseline and introduce steady-state reserve errors. On the other hand, a very large gain might introduce oscillations in the aggregate power and thus increased reserve errors. In fact, upper and lower bounds on $K_\text{c}$ can be computed analytically.

The correction of mean temperature with a gain $K_\text{c}$ can be thought as temperature resetting with the $K_{\text{r},t}$ defined in \eqref{resetting_factor2}, but for an opposite frequency deviation signal. Assuming that the application of the temperature resetting factor $K_{\text{r},t}$ does not create oscillations in aggregate power, and if the opposite temperature change due to $K_\text{c}$ is in absolute terms smaller than the average $\bar{K}_{\text{r},t}$, then no oscillations in aggregate power should occur. Therefore, an upper bound on $K_\text{c}$ can be obtained from
\begin{align}
&\big|K_\text{c} \cdot (\bar{T}_{t}-\bar{T}_{\text{nom}})\big| \leq \left|\bar{K}_{\text{r},t}\right| ~~\Leftrightarrow \nonumber\\
&K_\text{c} \hspace{-0.05cm}\cdot\hspace{-0.05cm} |\bar{T}_{t}\hspace{-0.05cm}-\hspace{-0.05cm}\bar{T}_{\text{nom}}| \hspace{-0.05cm}\leq\hspace{-0.05cm} \Delta{t} \bar{\beta} \bar{P}_\text{n} \bigg|\frac{D^{\textrm{r}}\bar{\Delta{f}}_t}{\Delta{f}_{\textrm{max}}}\bigg| \hspace{-0.05cm}=\hspace{-0.05cm} \Delta{t} \bar{\beta} \bar{P}_\text{n} \hspace{-0.05cm}\cdot\hspace{-0.05cm} |\bar{D}^\text{n}_t \hspace{-0.05cm}-\hspace{-0.05cm}D^\text{n}|.\label{Kc_and_duty_cycle}
\end{align}

From \eqref{duty_cycle_estimate1} we can see that $\bar{D}^\text{n}_t$ is a function of the mean temperature $\bar{T}_{t}$. Therefore, the upper bound on $K_\text{c}$ from \eqref{Kc_and_duty_cycle} is in general a function of $\bar{T}_{t}$ itself. However, a temperature-independent bound can be obtained by a first order Taylor expansion of \eqref{duty_cycle_estimate1} that gives
\begin{align}
&\bar{D}^\text{n}_t(\bar{T}) \approx \bar{D}^\text{n}_t\left(\bar{T}_{\text{nom}}\right) + \frac{d\bar{D}_t^\text{n}\left(\bar{T}\right)}{d\bar{T}}\Big|_{\bar{T}_\text{nom}} \cdot \left(\bar{T}-\bar{T}_\text{nom}\right) &~\Leftrightarrow \nonumber
\end{align}
\begin{align}
&\bar{D}^\text{n}_t\hspace{-0.06cm}-\hspace{-0.06cm}D^\text{n} \hspace{-0.06cm}=\hspace{-0.06cm} \bar{D}^\text{n}_t(\bar{T})\hspace{-0.06cm}-\hspace{-0.05cm}\bar{D}^\text{n}_t(\bar{T}_{\text{nom}}) \hspace{-0.06cm}\approx\hspace{-0.06cm} \frac{d\bar{D}_t^\text{n}\left(\bar{T}\right)}{d\bar{T}}\Big|_{\bar{T}_\text{nom}} \hspace{-0.45cm}\cdot\hspace{-0.06cm} (\bar{T}\hspace{-0.06cm}-\hspace{-0.06cm}\bar{T}_\text{nom}).\label{taylor_expansion}
\end{align}
Substituting now \eqref{taylor_expansion} in \eqref{Kc_and_duty_cycle} gives us the upper bound
\begin{align} \label{upper_Kc_bound}
K_\text{c} \leq \bigg|\Delta{t} \bar{\beta} \bar{P}_\text{n} \cdot \frac{d\bar{D}_t^\text{n}\left(\bar{T}\right)}{d\bar{T}}\Big|_{\bar{T}_\text{nom}}\bigg|~.
\end{align}

A lower bound on $K_\text{c}$ can be obtained based on the minimum/maximum acceptable mean temperature and the maximum acceptable settling time, i.e., the time needed to restore $\bar{T}_t$ close enough to $\bar{T}_{\text{nom}}$ after the biased frequency deviation event is over. Neglecting lockout constraints \eqref{mean_temp_est} can be written as
\begin{align}
\bar{T}_t &= \lambda \cdot \bar{T}_{t-1} - \gamma \cdot \Delta{f}_t + (1-\lambda) \cdot \bar{T}_{\text{nom}}~, \label{mean_temp_est_simple}
\end{align}
where $\lambda = 1-K_\text{c}$ and $\gamma = D^{\textrm{r}} \Delta{t} \bar{\beta} \bar{P_\text{n}}/\Delta{f}_{\textrm{max}}$. With an initial condition $\bar{T}_0$, \eqref{mean_temp_est_simple} has the solution
\begin{align}
\bar{T}_t &= \lambda^t \bar{T}_0 - \gamma \sum\nolimits_{k=0}^{t-1} \lambda^k \Delta{f}_{t-k-1} + (1-\lambda^t) \bar{T}_{\text{nom}}~. \label{mean_temp_est_simple_rec}
\end{align}

Assume that $\bar{T}_0 = \bar{T}_{\text{nom}}$ and that the refrigerator aggregation faces the step frequency deviation
\begin{align} \label{Df_signal}
\Delta{f}_t =
\begin{cases}
\delta, & \textrm{if} \hspace{2mm} 0 \leq t \leq N_\text{ev}\\
0, & \textrm{if} \hspace{2mm} t > N_\text{ev}
\end{cases}~,
\end{align}
which represents a frequency deviation signal with a bias equal to $\delta$ over a period $N_\text{ev}$. Let $\varepsilon$ denote the tolerance in terms of mean temperature during the frequency deviation event, i.e., $\bar{T}_t$ must satisfy $\bar{T}_{\text{nom}}-\varepsilon \leq \bar{T}_t \leq \bar{T}_{\text{nom}}+\varepsilon~\forall t$. Furthermore, let $N_\text{rec}$ denote the maximum acceptable settling time after the frequency deviation event, i.e, $\bar{T}_{\text{nom}}-\epsilon \leq \bar{T}_t \leq \bar{T}_{\text{nom}}+\epsilon$ should hold for $t\geq N_\text{ev}+N_\text{rec}$, where $\epsilon<\varepsilon$ is another tolerance. A lower bound on $K_\text{c}$ can be computed using Proposition~\ref{mean_temp_proposition}.

\begin{proposition} \label{mean_temp_proposition}
If $\bar{T}_0 = \bar{T}_{\text{nom}}$ and $\Delta{f}_t$ is given by \eqref{Df_signal}, then $\bar{T}_{\text{nom}}-\varepsilon \leq \bar{T}_t \leq \bar{T}_{\text{nom}}+\varepsilon~\forall t$ and $\bar{T}_{\text{nom}}-\epsilon \leq \bar{T}_t \leq \bar{T}_{\text{nom}}+\epsilon,~\text{for } t\geq N_\text{ev}+N_\text{rec},$ hold if the corrective temperature gain is computed as $K_\text{c}=1-\lambda$, where $\lambda$ satisfies
\begin{align}
\gamma \delta \cdot (1-\lambda^{N_\text{ev}}) &\leq \varepsilon \cdot (1-\lambda) \label{lambda_ineq1} \\
\gamma \delta \cdot \lambda^{N_\text{rec}} \cdot (1-\lambda^{N_\text{ev}}) &\leq \epsilon \cdot (1-\lambda) \label{lambda_ineq2}~.
\end{align}
\end{proposition}

\begin{proof}
The proof is given in the appendix.
\end{proof}

\begin{remark}\label{remark_mean_temp}
Although \eqref{lambda_ineq1}, \eqref{lambda_ineq2} cannot be solved analytically, a lower bound on $K_\text{c}$ is obtained numerically starting with $\lambda=1$ and gradually reducing it until \eqref{lambda_ineq1}, \eqref{lambda_ineq2} are satisfied.
\end{remark}

Equation~\eqref{upper_Kc_bound}, Proposition~\ref{mean_temp_proposition}, and Remark~\ref{remark_mean_temp} provide theoretical lower and upper bounds on $K_\text{c}$ that can serve as an initial range when determining the final value of the gain via simulations with realistic $\Delta{f}_t$ signals, and in the presence of startup dynamics and lockout constraints.

\subsection{Addressing Limited Thermostat Resolution}\label{thermostat_theory}
Typically, the change in the thermostat's temperature limits from \eqref{Tmin_corrective} and \eqref{Tmax_corrective} is very small, for example, $10^{-4}~^\circ\textrm{C/s}$. Nevertheless, the measurement accuracy in a refrigerator can be as low as $\pm 0.2^\circ\textrm{C}$ \cite{sossan2016grey}. Therefore, the temperature limits of some refrigerators might not change, and thus the aggregate power might decay towards its pre-disturbance value. 

We address this practical issue by allowing the temperature limits to change according to a sufficiently large fixed rate $\Delta{T}_\text{res}~^\circ\textrm{C}$ per second. Specifically, instead of distributing the temperature change over the whole population and requesting a small change from every device, we impose a larger change on the temperature limits of fewer devices.

\subsubsection{Probabilistic Implementation}
The required change in temperature limits is
\begin{align}
\Delta{T}_{\text{lim},t} = K_{\text{r},t} - K_\text{c} \cdot (\bar{T}_{t-1}-\bar{T}_{\text{nom}})~.
\end{align}
If the temperature limits of an individual device change only by $\Delta{T}_\text{res}$, we can get the same total change in temperature by applying $\Delta{T}_\text{res}$ to a fraction of the population equal to
\begin{align}\label{prob_temp_resetting}
\varrho_t = |\Delta{T}_{\text{lim},t}|/\Delta{T}_\text{res}~.
\end{align}
The implementation of this approach is straightforward: each refrigerator draws a random number uniformly distributed between $0$ and $1$ and if it is smaller than $\varrho_t$ and the refrigerator is unlocked, then it modifies its temperature limits by $\Delta{T}_\text{res}$, otherwise it does not. The limits increase by $\Delta{T}_\text{res}$ if $\Delta{T}_{\text{lim},t}>0$ and decrease by $\Delta{T}_\text{res}$ if $\Delta{T}_{\text{lim},t}<0$.

Although simple, this approach introduces two sources of inaccuracy in the control. First, as $\Delta{T}_\text{res}$ increases, fewer devices change their thermostat limits and therefore the control is exposed to inaccuracies from random number generation. Second, large $\Delta{T}_\text{res}$ values might synchronize the population. Without the fixed rate $\Delta{T}_\text{res}$, the temperature limits change in a continuous way, the devices switch when they are close to the thermostat limits and this avoids synchronization. In contrast, if a large $\Delta{T}_\text{res}$ value is used, some devices might switch at a temperature considerably different than the thermostat limits.

\subsubsection{Effect on Temperature Limits and Countermeasure}
If a fixed rate $\Delta{T}_\text{res}$ is used, the change in temperature limits $\Delta{T}_{\text{lim},i,t}$ is different for each refrigerator $i$. A disadvantage of this approach is that the variance of $\Delta{T}_{\text{lim},i,t}$ monotonically increases, as shown by Proposition~\ref{var_temp_proposition} and Remark~\ref{var_temp_remark}.

\begin{proposition} \label{var_temp_proposition}
If the refrigerators' temperature limits change probabilistically according to \eqref{prob_temp_resetting}, then the mean and variance of $\Delta{T}_{\text{lim},i,t}$ are given by
\begin{align}
\Mean\left[\Delta{T}_{\text{lim},i,t}\right] &= \sum\nolimits_{k=0}^{t-1} \Delta{T}_{\text{lim},k} \label{mean_with_DTres}\\
\mathrm{Var}\left[\Delta{T}_{\text{lim},i,t}\right] \hspace{-0.06cm}&=\hspace{-0.06cm} \sum\nolimits_{k=0}^{t-1} |\Delta{T}_{\text{lim},k}| \hspace{-0.04cm}\cdot\hspace{-0.04cm} \left(\Delta{T}_\text{res}\hspace{-0.04cm}-\hspace{-0.04cm}|\Delta{T}_{\text{lim},k}|\right). \label{var_with_DTres}
\end{align}
\end{proposition}

\begin{proof}
The proof is given in the appendix.
\end{proof}

\begin{remark}\label{var_temp_remark}
If no fixed rate $\Delta{T}_\text{res}$ is applied, all refrigerators change their temperature limits by $\sum_{k=0}^{t-1} \Delta{T}_{\text{lim},k}$, so the mean value is $\sum_{k=0}^{t-1} \Delta{T}_{\text{lim},k}$ and the variance is zero. In practice, the mean value will be slightly different and the variance non-zero but small because some devices are always locked at steady-state and do not change their temperature limits. According to \eqref{mean_with_DTres}, $\Delta{T}_\text{res}$ does not affect the mean value of temperature limit change. In contrast, the variance monotonically increases in time (linearly with $\Delta{T}_\text{res}$) because $\Delta{T}_\text{res}-|\Delta{T}_{\text{lim},k}| \geq 0$.
\end{remark}

Due to the monotonic increase in variance, some devices will likely sustain larger and prolonged deviations from the nominal temperature limits, which is undesirable. The variance can be bounded by imposing a bound $\Delta{T}_{\text{b}}$ on the maximum deviation from the mean temperature of the population. To implement this, we require each unlocked device $i$ that drew a random number smaller than $\varrho_t$ at time step $t$ to change its limits only if
\begin{align} \label{var_bound}
\bar{T}_{t-1}-\Delta{T}_{\text{b}} \leq \Delta{T}_{\text{lim},i,t-1} + \zeta_i \Delta{T}_\text{res} \leq \bar{T}_{t-1}+\Delta{T}_{\text{b}}~,
\end{align}
where $\zeta_i$ is equal to $1$ or $-1$ depending on the sign of $\Delta{T}_{\text{lim},t}$.

\section{Refrigerator Door Openings} \label{door_openings}
The analysis so far neglected the door openings that significantly affect the aggregate power of a population of refrigerators. Since the PFC performance is evaluated with respect to a baseline, neglecting the effect of door openings will introduce large control errors. In this section, we present a method to account for door openings in \ac{PFC} design.

\subsection{Stochastic Model for Door Opening Events}
There exist only a few papers that consider door openings in power system studies with refrigerators, e.g., \cite{Keep2011switches,kremers2013emergent}. We assume that the number of door openings per day follows a normal distribution with mean value $\mu_\text{op}$ and standard deviation $\sigma_\text{op}$. Furthermore, we assume that the duration of a door opening event follows a normal distribution with mean value $\mu_\text{d}$ and standard deviation $\sigma_\text{d}$. Based on \cite{parker1992openings} we select $\mu_\text{op}=40$, whereas we fix $\sigma_\text{op}=5$, $\mu_\text{d}=20$~s and $\sigma_\text{d}=3$~s \cite{kremers2013emergent}.

\subsection{Modeling the Effect of Door Openings}
The door openings increase the refrigerator's energy consumption by $\xi \cdot 100\%$ (we assume $\xi=0.22$ \cite{parker1992openings}). The effect of door openings on refrigerator temperature is modeled by reducing the thermal resistance of the model \eqref{eq:ode_off}, \eqref{eq:ode_on} during each door opening event, i.e., by increasing the thermal losses to the ambient. An upper bound to the new thermal resistance is given by Proposition~\ref{door_opening_prop}.

\begin{proposition}\label{door_opening_prop}
Denote by $N^\text{d}=86,400$~s the duration of a day. An upper bound to the estimate of the thermal resistance $R_\text{op}$ during a door opening event can be obtained from
\begin{align}\label{Ropen}
R_\text{op} \leq R \cdot \frac{1}{1+\frac{N^\text{d}}{\mu_\text{op} \cdot \mu_\text{d}} \cdot \xi}~.
\end{align}
\end{proposition}

\begin{proof}
The proof is given in the appendix.
\end{proof}

\begin{remark}
Evaluating \eqref{Ropen} as equality gives us an initial estimate of $R_\text{op}$. For the assumed parameters, \eqref{Ropen} gives \mbox{$R_{\text{op}} \leq R/24.76$}, and so we set $R_{\text{op}}=R/25$ in our simulations.
\end{remark}

\subsection{Modifications in Primary Frequency Control Design}
The proposed controller can be easily extended to account for door openings. Even without \ac{PFC}, the aggregation's duty cycle and baseline power are not constant any more, but they depend on the distribution of door openings within the day. The duty cycle without \ac{PFC} can be expressed as
\begin{align}
D^\text{n}_t = D^\text{n} + \Delta{D}^\text{op}_t~,
\end{align}
where $\Delta{D}^\text{op}_t\geq0$ is the additional duty cycle due to door openings. A smoothed version of $\Delta{D}^\text{op}_t$ can be obtained comparing historical data of aggregate power consumption with and without door openings.

Instead of \eqref{duty_cycle_estimate1}, the duty cycle is now estimated using
\begin{align}
\bar{D}^\text{n}_t \approx \frac{\mathrm{E}\left[t_{\text{on},t}\right]}{\mathrm{E}\left[t_{\text{on},t}\right]+\mathrm{E}\left[t_{\text{off},t}\right]} + \Delta{D}^\text{op}_t~.
\end{align}
The necessary information about door openings is passed to the switching probability calculation \eqref{eq:rho} through \eqref{XandD2}, and no other modification in the control design is needed. In particular, note that each refrigerator responds to frequency deviations irrespective of whether the door is closed or open.

\begin{figure*}[t]
\centering
\includegraphics[width=1.2\textwidth,height=2in]{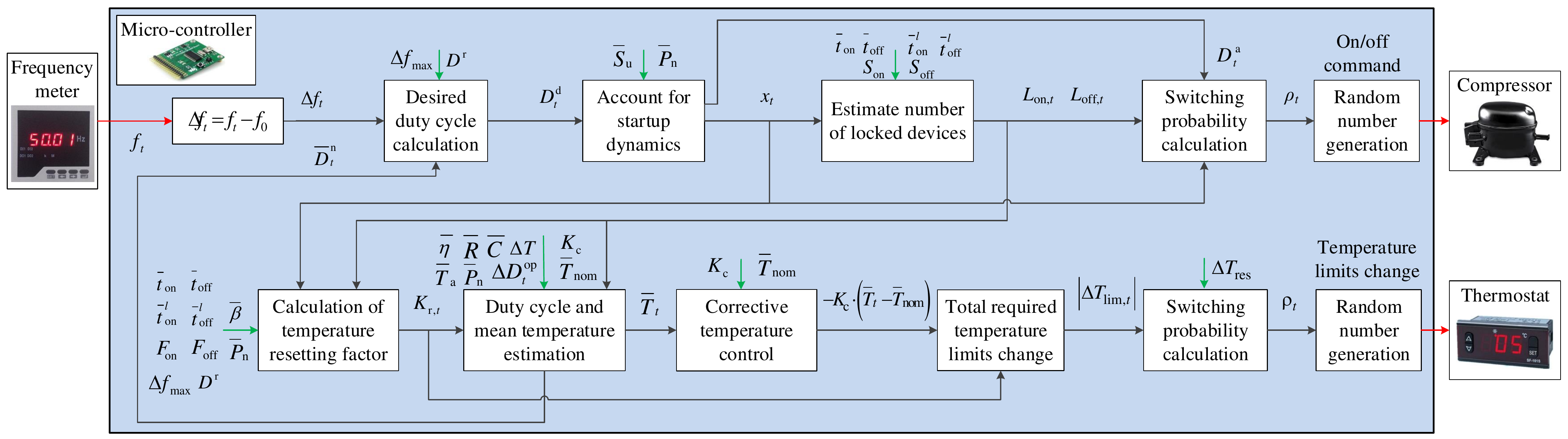}
\caption{The proposed decentralized stochastic controller. Red arrows indicate the controller's input/output signals, green arrows indicate input from stored parameter values in the controller's memory, and black arrows indicate internal control and feedback signals.} \label{fig:control_overview}
\end{figure*}

\section{Parameters and Performance Metrics} \label{parameters_and_metrics}
The proposed controller (consisting of all components presented in the previous sections) is graphically shown in Fig.~\ref{fig:control_overview} and is easy to implement at the device level. The required hardware includes a frequency meter and a micro-controller to control the switch of the compressor and the temperature limits of the thermostat. The micro-controller consists of a random number generator and a floating-point unit able to perform elementary mathematical operations (additions, subtractions, multiplications and divisions).

In this section, we present the refrigerator parameters in Table~\ref{table:general_fridge_params}, and introduce the metrics and benchmarks that we use to quantify the controller's performance. Note that $\beta$ and $t_{\textrm{off}}^{l}$ have the mean value and standard deviation shown in Table~\ref{table:general_fridge_params}, but they do not follow a standard \ac{PDF} because they are post-calculated based on the imposed \acp{PDF} on other variables \cite{Ziras2015GM}.

\begin{table}[t]
\renewcommand{\arraystretch}{1.15}
\caption{Refrigeration Aggregation Parameters}
\centering
\begin{tabular}{c c | c c}
\hline
Parameter & Value & Parameter & Value\\%[0.1ex]
\hline\
$T_{\textrm{a}}$~(\textrm{$^\circ$}C) & $\Unif[20, 24]$ & $\alpha~(\textrm{s}^{-1})$ & $\Unif[4, 6]\cdot10^{-5}$\\
$\Delta{T}$~(\textrm{$^\circ$}C) & $\Unif[1.7, 2.3]$ & $\beta$~(\textrm{$^\circ$}C/\textrm{kJ}) & $\left(\mu, \sigma\right)=(4.4, 0.7)\cdot10^{-5}$\\
$T_\text{set}$~(\textrm{$^\circ$}C) & $\Unif[4.5, 5.5]$ & $P_{\textrm{n}}$~(W) & $\Unif[70, 90]$\\
$u$~(-) & $\Normal[0.25, 0.025]$ & $t_{\textrm{on}}^{l}$~(s) & $\Normal[60, 5]$\\
$N_\text{s}~(s)$ & $\Normal[30, 3]$ & $t_{\textrm{off}}^{l}$ (s) & $(\mu, \sigma)=(189, 31.5)$\\
\hline
\end{tabular}
\label{table:general_fridge_params}
\end{table}

For a fixed $D^\text{r}$ the PFC reserve capacity is given by
\begin{align}
P_{\text{res}} = N_\textrm{r} \bar{P}_{\textrm{n}} D^\text{r}~.
\end{align}
The desired aggregate power at time step $t$ with PFC is
\begin{align}
P_{\text{d},t} = \bar{P}_{\text{b},t} + P_{\text{res}} \cdot (\Delta{f}_t/\Delta{f}_{\textrm{max}})~,
\end{align}
where $\bar{P}_{\text{b},t}$ is the population's baseline power, i.e., a smoothed version of historical measurements of aggregate power $P_{\text{b},t}$ without PFC. If $w_{i}(t)=0$~$\forall i$, the baseline is constant and equal to $\bar{P}_\text{b}=N_\textrm{r} \bar{P}_{\textrm{n}} D^\text{n}$.

The reserve \ac{MAPE} $e_\text{r,mape}$ and the tracking MAPE $e_\text{t,mape}$ are used as control performance metrics. Note that in the definition of tracking MAPE the normalized error is computed by dividing the absolute error with the desired power, whereas in the reserve MAPE definition the normalization is performed using the reserve capacity. In addition, we use the baseline \ac{MAPE} $e_\text{b,mape}$ as a metric of the natural oscillations in an uncontrolled refrigerator aggregation. The metrics are defined as ($N_\text{sim}$ is the simulation period)
\begin{align}
e_\text{r,mape} &= (100/N_\text{sim}) \cdot \sum\nolimits_{t=0}^{N_\text{sim}-1} \left|(P_{\text{d},t}-P_{\text{agg},t})/P_{\text{res}}\right| \label{res_mape}
\end{align}
\begin{align}
e_\text{t,mape} &= (100/N_\text{sim}) \cdot \sum\nolimits_{i=0}^{N_\text{sim}-1} \left|(P_{\text{d},t}-P_{\text{agg},t})/P_{\text{d},t}\right|  \label{track_mape}
\end{align}
\begin{align}
e_\text{b,mape} &= (100/N_\text{sim}) \cdot \sum\nolimits_{t=0}^{N_\text{sim}-1} \left|(\bar{P}_{\text{b},t}-P_{\text{b},t})/P_{\text{res}}\right|.\label{baseline_mape}
\end{align}

The \emph{``proposed controller''} is benchmarked against a \emph{``simple controller 1''} that neglects the startup dynamics and lockout constraints, and is similar to the approach of \cite{borsche2015new}. The only difference is that the temperature resetting factor $K_{\text{r},t}$ is computed analytically from \eqref{resetting_factor2} and not chosen heuristically as in \cite{borsche2015new}. We also consider a \emph{``simple controller 2''} that apart from neglecting the startup dynamics and lockouts, it does not reset the thermostat limits as well.

The controllers' comparison is performed via detailed simulations using actual frequency deviation data from Switzerland in $2009$ and in $2011$. We assume that the refrigerators participate with a small share in PFC, and therefore have a negligible effect on system frequency. For this reason, it is sufficient to model the frequency as a time series in our simulations in order to evaluate the controller's performance. Of course, a dynamic frequency model incorporating also the generator's power response is needed for high shares of refrigerators in PFC, because their effect on frequency cannot be neglected anymore. In \cite{Vrettos2015PowerTechConf} we present results from dynamic frequency studies to show the contribution of refrigerators in reducing under-frequency after sudden loss of generation.

\section{Sensitivity to Design and Load Parameters} \label{sensitivity}
\subsection{Sensitivity to Aggregation Size} \label{aggregation_size_sims}
The control is expected to perform better for large refrigerator aggregations for three main reasons. First, the larger the aggregation the closer the fraction of switched loads to the desired switching probability. Second, the controller is designed based on average refrigerator parameters. And third, the assumption of a smooth uncontrolled duty cycle $D^\text{n}_t$ is reasonable only for sufficiently large aggregations, where the oscillations from individual load cycles cancel out.

\begin{figure}[t]
\centering
\includegraphics[width=\textwidth]{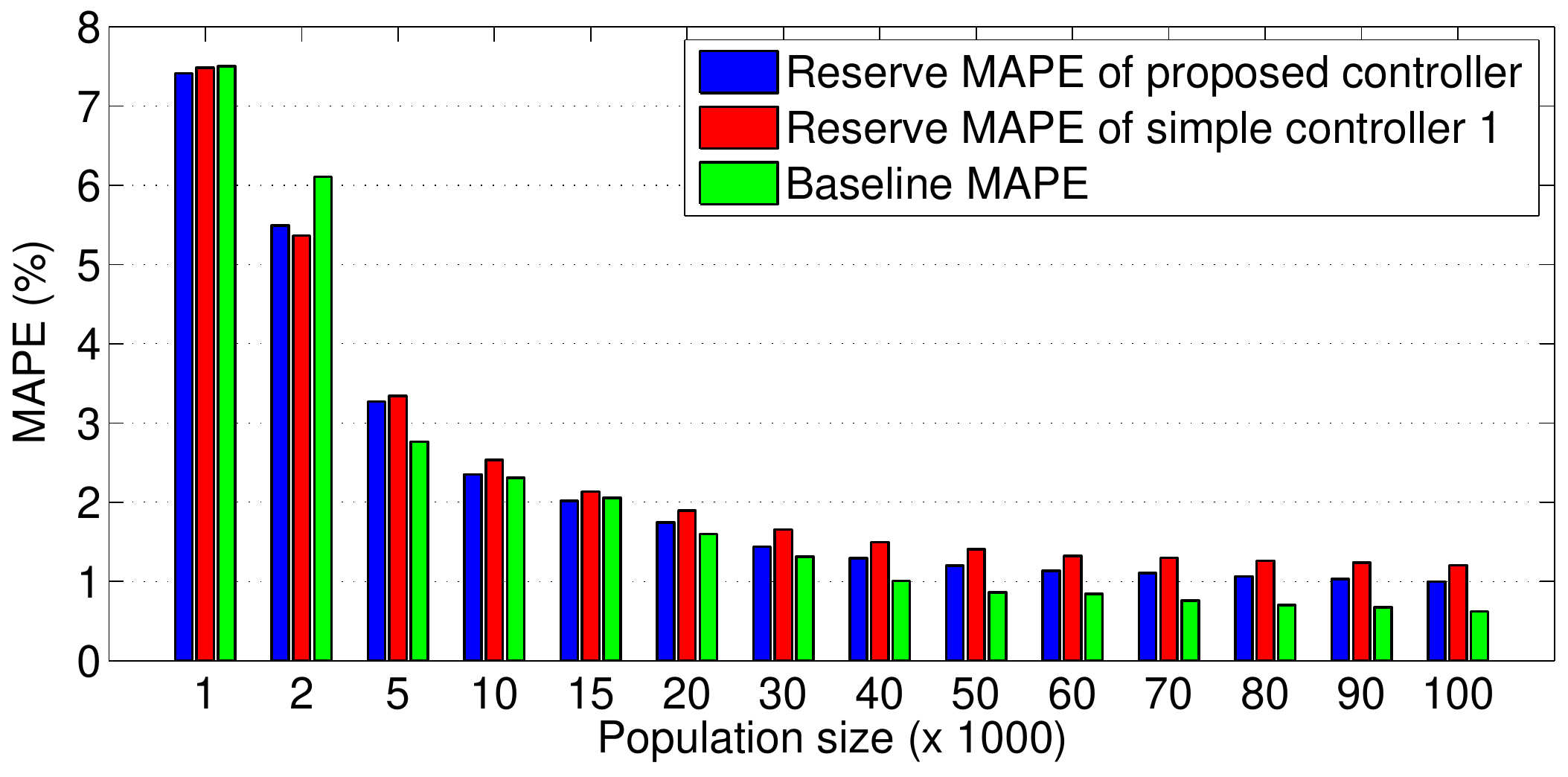}
\caption{Dependence of control performance and baseline error on the size of refrigerator aggregation for the proposed controller and simple controller 1.} \label{fig:reserve_error_vs_aggr_size}
\end{figure}

We investigate the dependence of controller performance on aggregation size by running simulations using $30$ zero-mean, $5$-hour samples of frequency deviation ($D^\text{r}=0.15$ is used). Figure~\ref{fig:reserve_error_vs_aggr_size} shows the reserve MAPE ($e_\text{r,mape}$) for the proposed and the simple controller 1, which decrease with the square root of the aggregation size. A population of $1,000$ refrigerators results in a large reserve error, which drastically decreases when the size increases to $10,000$. After this point, the reserve error decreases asymptotically and it practically saturates at an aggregation size of $70,000$. The proposed controller outperforms the simple controller 1 for all aggregation sizes apart from an aggregation of $2,000$ refrigerators. For populations of $70,000$ loads or more, the proposed controller reduces the reserve error by $15\%$.

Simple controller 1 might perform better than the proposed controller for small aggregations, due to the dominant effect of the baseline MAPE ($e_\text{b,mape}$) for aggregations up to $15,000$ loads, as shown in Fig.~\ref{fig:reserve_error_vs_aggr_size}. In other words, the PFC response of small aggregations is covered by the large natural oscillations of the uncontrolled baseline. Therefore, the improvement of the proposed controller over simple controller 1 is guaranteed for sufficiently large aggregations, but not for small ones. Interestingly, the baseline error comprises a significant part of the control error even for large aggregations (for example, approximately $62\%$ of the total error for $100,000$ loads).

\subsection{Sensitivity to Reserve Capacity}
We investigate the controller's sensitivity to reserve capacity $D^\text{r}$ for an aggregation of $70,000$ refrigerators (according to the results of Fig.~\ref{fig:reserve_error_vs_aggr_size}). We ran simulations for $10$ days with different values of $D^\text{r}$ using as base case scenario the parameters of Table~\ref{table:general_fridge_params}. The frequency deviation is approximately zero-mean during the first $5$ days of the data set, whereas it has a significant bias during the remaining $5$ days. The same frequency data are used for the all analyses of the rest of Section~\ref{sensitivity}.

\begin{figure}[t]
\centering
\includegraphics[width=\textwidth]{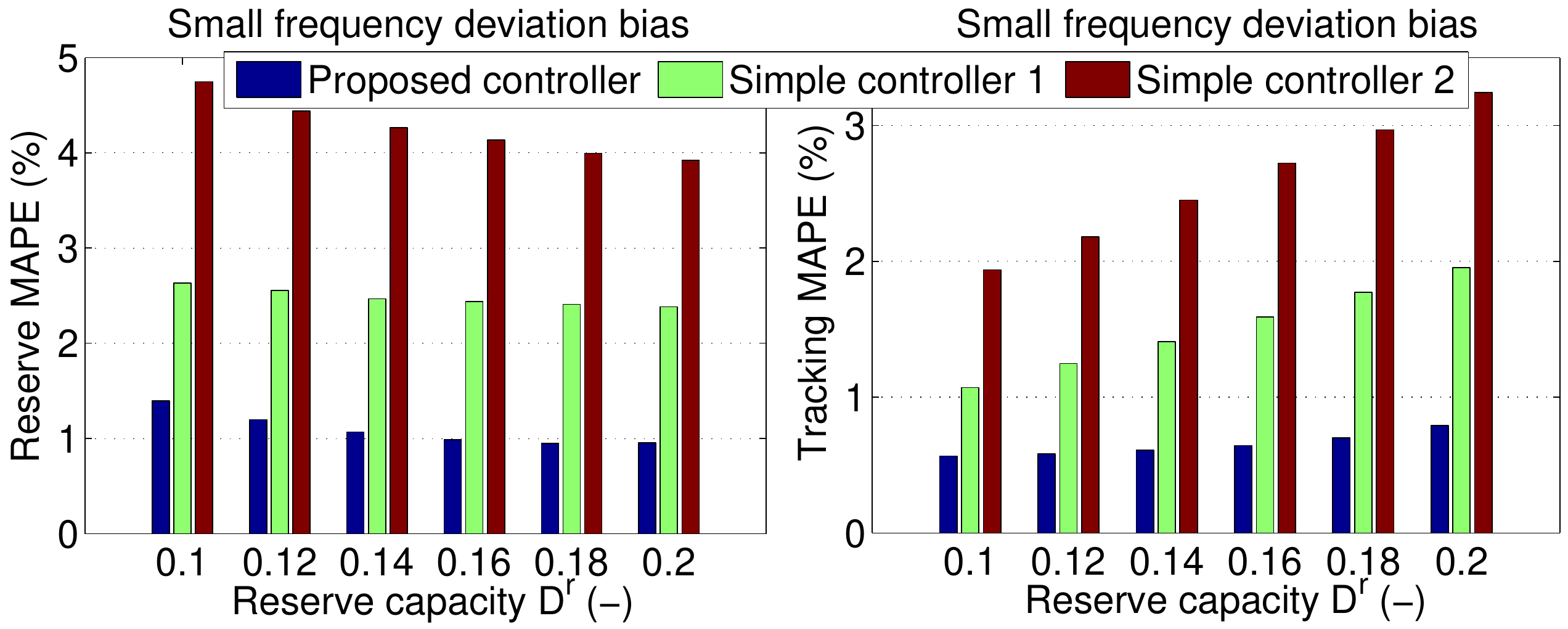}
\caption{Sensitivity to reserve capacity in simulations with small frequency deviation bias. Left: Reserve MAPE. Right: Tracking MAPE.} \label{fig:capacity_results_5days_small}
\end{figure}

\begin{figure}[t]
\centering
\includegraphics[width=\textwidth]{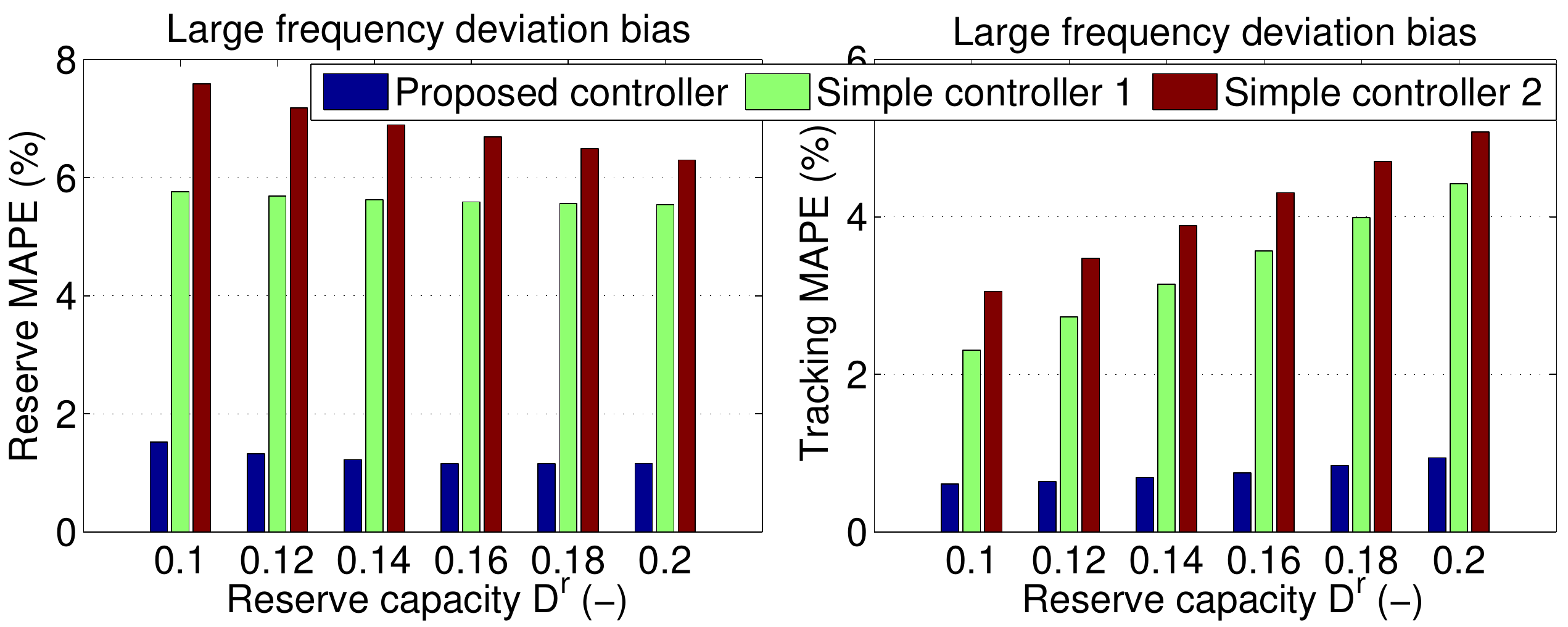}
\caption{Sensitivity to reserve capacity in simulations with large frequency deviation bias. Left: Reserve MAPE. Right: Tracking MAPE.} \label{fig:capacity_results_5days_large}
\end{figure}

Figure~\ref{fig:capacity_results_5days_small} shows the sensitivity analysis results for the frequency deviations with small bias, whereas Fig.~\ref{fig:capacity_results_5days_large} is for the frequency deviations with large bias ($K_\text{c}=0.5\cdot 10^{-4}$ is used based on the controller tuning results of Section~\ref{empirical_Kc}). In both figures, the left plot shows the reserve MAPE and the right plot the tracking MAPE.

The reserve capacity that minimizes the control error depends on the definition of the error metric. Specifically, as $D^\text{r}$ increases the tracking MAPE increases but the reserve MAPE decreases. Therefore, although increasing the reserve capacity generally results in a higher control error, the error becomes smaller in relative terms (i.e., in comparison with the reserve capacity that is offered). The proposed controller consistently reduces the errors compared with simple controller $1$. In addition, the latter performs better than simple controller $2$, which illustrates the need for resetting of thermostat temperature limits.

Observe that the difference between the proposed controller and simple controller $1$ is more pronounced for biased frequency deviations. This shows the effectiveness of the corrective temperature control of Section~\ref{corrective_control_biased_df}. In fact, there is only a small difference between the simple controller $1$ and simple controller $2$ for biased frequency deviations, because in both cases most of the reserve error is due to steady-state baseline errors.

\subsection{Sensitivity to Peak Power of Startup Dynamics}
\begin{figure}[t]
\centering
\includegraphics[width=\textwidth]{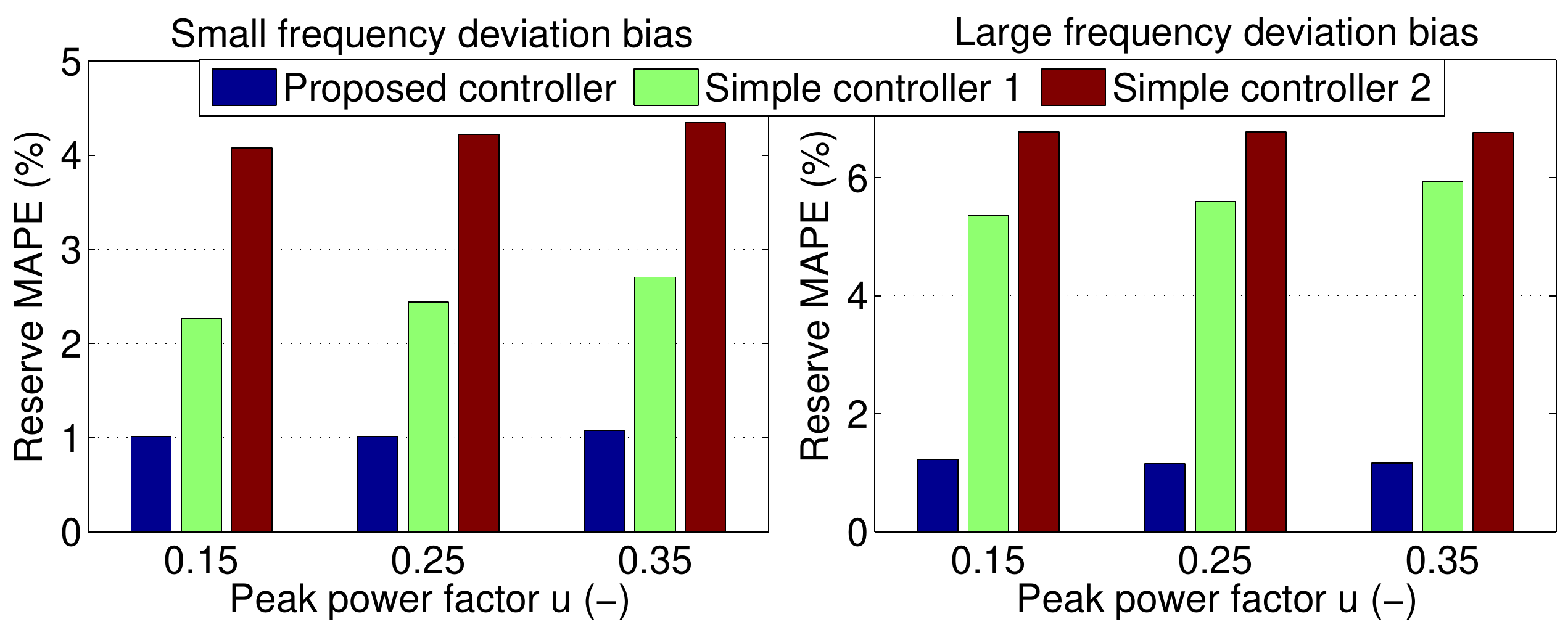}
\caption{Sensitivity of reserve MAPE to peak power during startup dynamics. Left: Small frequency deviation bias. Right: Large frequency deviation bias.} \label{fig:dynamics_results_10days}
\end{figure}

The dependence of reserve MAPE on the peak power factor ($u$) during the startup phase is presented in Fig.~\ref{fig:dynamics_results_10days}. The error of the proposed controller is not sensitive to $u$, whereas the error of the simple controller $1$ clearly increases as $u$ increases. This is expected because larger $u$ values result in larger overshoots in aggregate power consumption, if the effect of startup dynamics is not considered in the control design.

\subsection{Sensitivity to Lockout Times}
\begin{figure}[t]
\centering
\includegraphics[width=\textwidth]{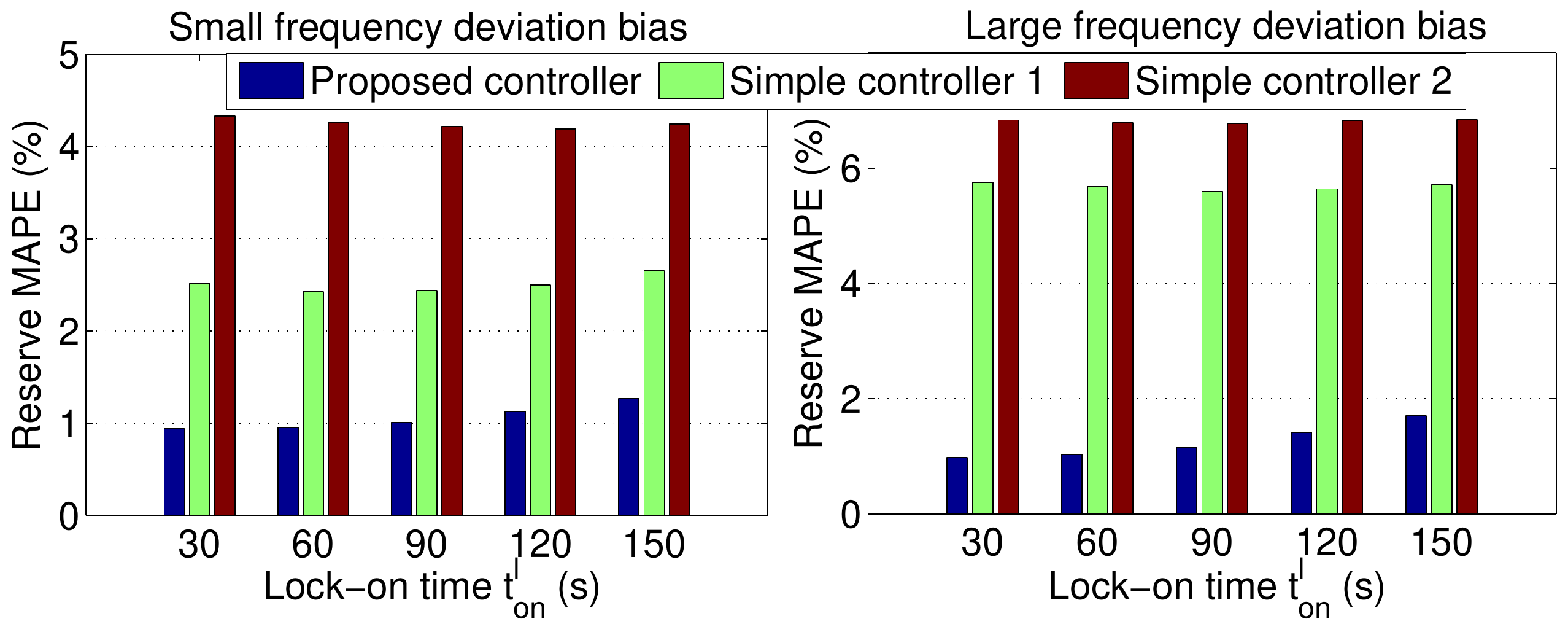}
\caption{Sensitivity of reserve MAPE to lock-on time. Left: Small frequency deviation bias. Right: Large frequency deviation bias.} \label{fig:lockout_results_10days}
\end{figure}

Figure~\ref{fig:lockout_results_10days} shows the dependence of reserve error to the lock-on time $t_\text{on}^l$. The error of the proposed controller steadily increases as $t_\text{on}^l$ increases, which indicates that the estimates of the number of locked devices -- obtained by \eqref{eq:countONOFF} -- deteriorate. The error of simple controller $1$ is significantly higher than that of the proposed controller, but is practically unaffected by $t_\text{on}^l$. This happens because simple controller $1$ neglects both the lockouts and startup dynamics. With reference to \eqref{eq:rho}, simple controller $1$ activates fewer loads for PFC than actually needed, because it neglects the lockouts. However, each of the activated loads consumes more power than assumed by the controller due to startup dynamics, which counteracts the error in the number of activated loads.

\subsection{Effect of PFC Activation Deadband}
In practice, a frequency deadband is used to prevent PFC activation to the wrong direction due to frequency measurement error. To investigate the effect of frequency deadband on the performance of the proposed controller, we repeated the simulations for our base case scenario with the deadband $\pm 10$~mHz of the ENTSO-E system \cite{swissgridPFC}.

\begin{figure}[t]
\centering
\begin{minipage}{0.49\linewidth}%
\centering \includegraphics[width=0.98\textwidth]{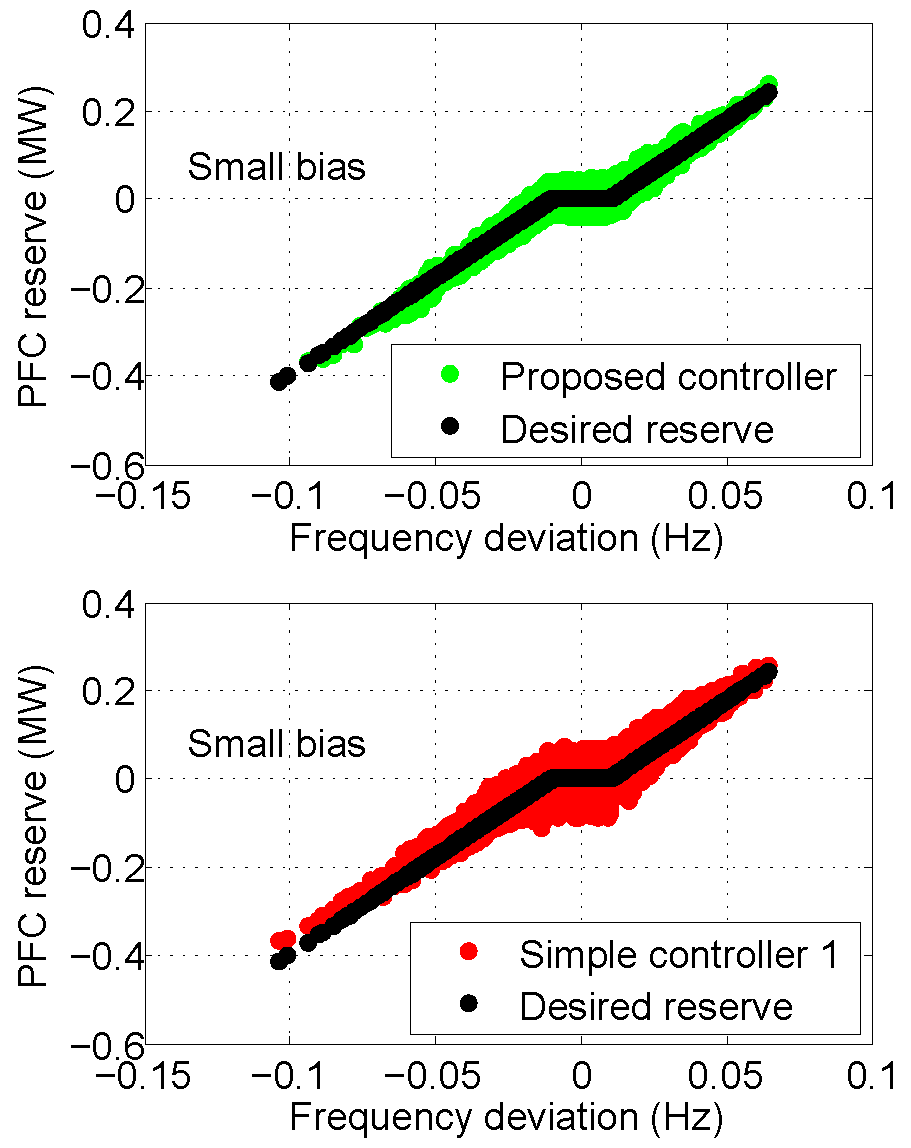}
\end{minipage}
\begin{minipage}{0.49\linewidth}
\centering \includegraphics[width=0.98\textwidth]{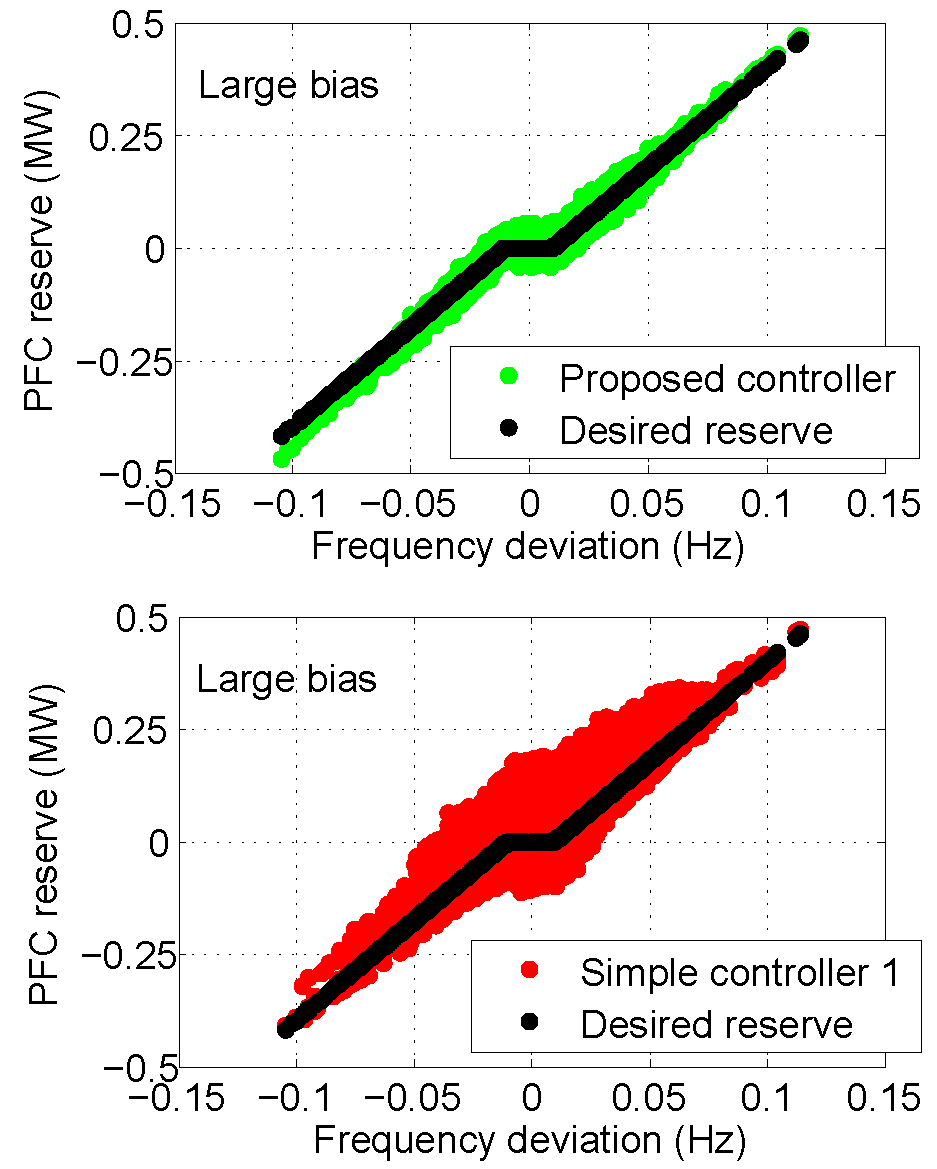}
\end{minipage}
\caption{Droop characteristic with a PFC activation deadband. Left: Small frequency deviation bias. Right: Large frequency deviation bias.}
\label{fig:droop_characteristic}
\end{figure}

Even with a deadband, the proposed controller significantly reduces the reserve error compared with simple controller $1$; however, the difference between the controllers is smaller because the PFC is active only for approximately $60\%$ of the time. Figure~\ref{fig:droop_characteristic} shows the desired and activated reserve for the two controllers with a PFC deadband for small and large frequency deviation bias. The proposed controller follows the desired droop characteristic of PFC much closer than simple controller $1$, in particular for biased frequency deviations where the error is reduced from $4.3\%$ to $1.2\%$.

\subsection{Evaluation}
Table~\ref{table:average_MAPE} summarizes the average values of reserve MAPE computed across all simulation sets. The reduction in control error achieved by the proposed controller (in comparison with the simple controllers) depends strongly on the frequency deviation bias. If the frequency deviation is zero-mean, the error reduces by approximately $15\%$. For small frequency deviation bias, the proposed controller reduces the error by approximately $56\%$ compared with simple controller $1$ and $74\%$ compared with simple controller $2$. For large frequency deviation bias, the error reduction is approximately $78\%$ compared with simple controller $1$ and $82\%$ compared with simple controller $2$.

\begin{table}[t]
\renewcommand{\arraystretch}{1.05}
\caption{Average Reserve MAPE ($\%$) for Different Frequency Biases}
\centering
\begin{tabular}{c|c|c|c}
\hline
Controller & Zero-mean & Small bias & Large bias\\
\hline
A. Simple controller $2$ & $1.34$ & $4.21$ & $6.76$ \\
B. Simple controller $1$ & $1.30$ & $2.48$ & $5.55$ \\
C. Proposed controller & $1.11$ & $1.08$ & $1.24$ \\
\hline
Improvement ($\%$) B $\rightarrow$ C & $14.62$ & $56.45$ & $77.66$ \\
Improvement ($\%$) A $\rightarrow$ C & $17.16$ & $74.35$ & $81.66$ \\
\hline
\end{tabular}
\label{table:average_MAPE}
\end{table}

The reserve MAPE of simple controller $1$ is already quite low in absolute terms, if the frequency deviation bias is small. Even in this case, however, the improvement achieved by the proposed controller is valuable in order to pass the pre-qualification tests that many system operators apply for reserve providers. In addition, the proposed controller will increase the profits from offering high-quality PFC reserves in a reserve market with performance-based payments.

Figure~\ref{fig:time_domain_plot} presents time-domain simulation results to compare the proposed controller against simple controller $1$. Specifically, we present the time series of frequency deviation, refrigerator aggregation's power response, and aggregation's mean temperature. As shown in Figs.~\ref{fig:time_domain_plot}b and~\ref{fig:time_domain_plot}d, there are notable offsets (long-term accumulated errors) in the response of simple controller $1$, but not in that of the proposed controller. This significant improvement in ``continuous operation'' of PFC with refrigerators is due to our contributions in control design. On the other hand, as shown in Fig.~\ref{fig:time_domain_plot}e, the instantaneous ``dynamic response'' of the proposed controller is only marginally better than that of simple controller $1$. Moreover, the proposed controller improves user utility by maintaining the aggregation's mean temperature close to the nominal value of $5^\circ$C (see Fig.~\ref{fig:time_domain_plot}c). This is very important in practice to ensure the participants' commitment to the demand response program.

\begin{figure}[t]
\centering
\includegraphics[width=0.77\textwidth,height=4.7in]{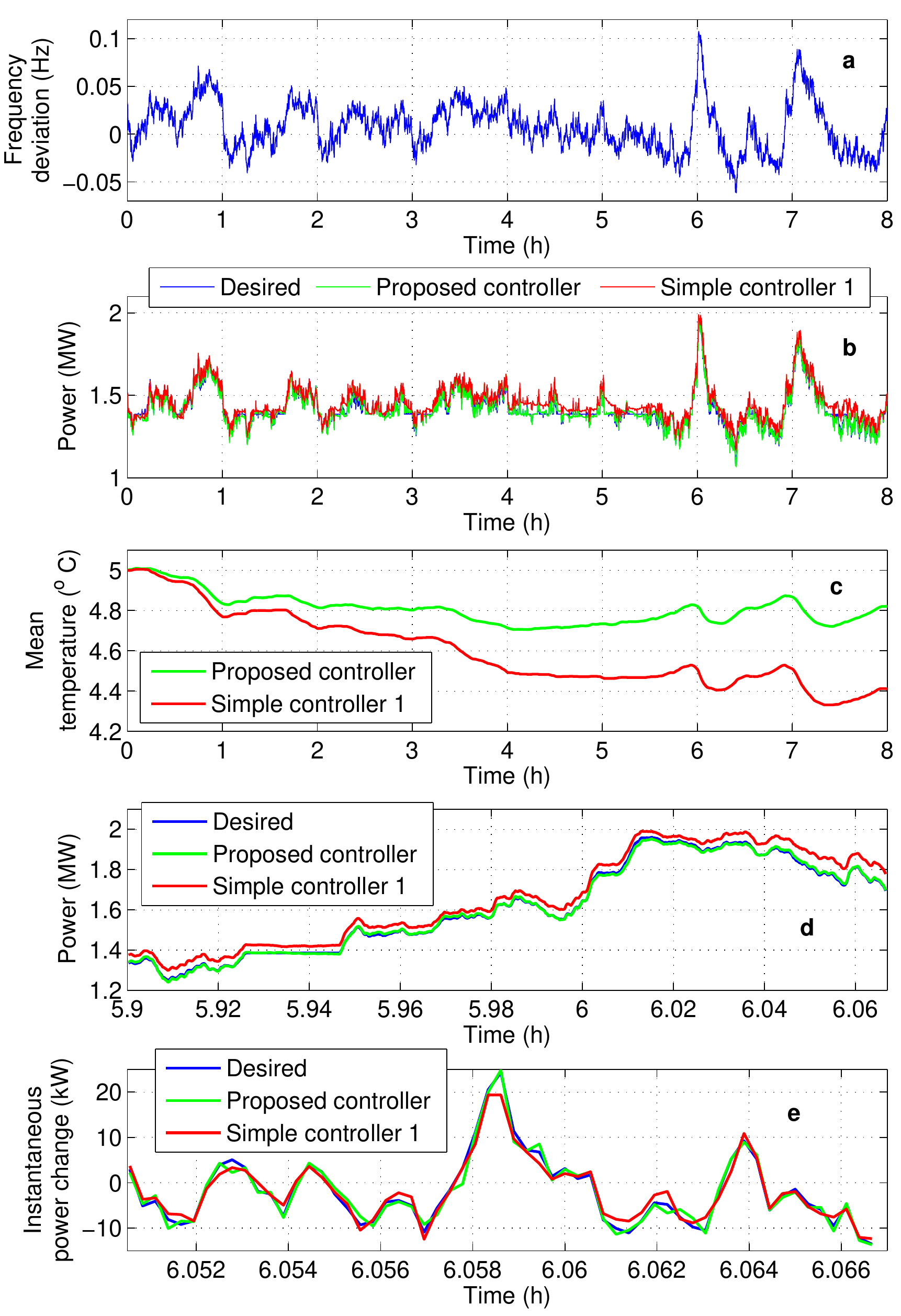}
\caption{Comparison of the proposed controller against simple controller $1$: (a) Frequency deviation; (b) Aggregate refrigerator power; (c) Mean temperature in the aggregation; (d) 10-minute close-up of plot (b); (e) 1-minute close-up of plot (d) with the instantaneous power change to provide PFC reserve.} \label{fig:time_domain_plot}
\label{fig:time_domain_plot}
\end{figure}

\section{Robustness Analysis} \label{robustness_analysis}
\subsection{Performance with Biased Frequency} \label{empirical_Kc}
\begin{figure}[t]
\centering
\includegraphics[width=\textwidth]{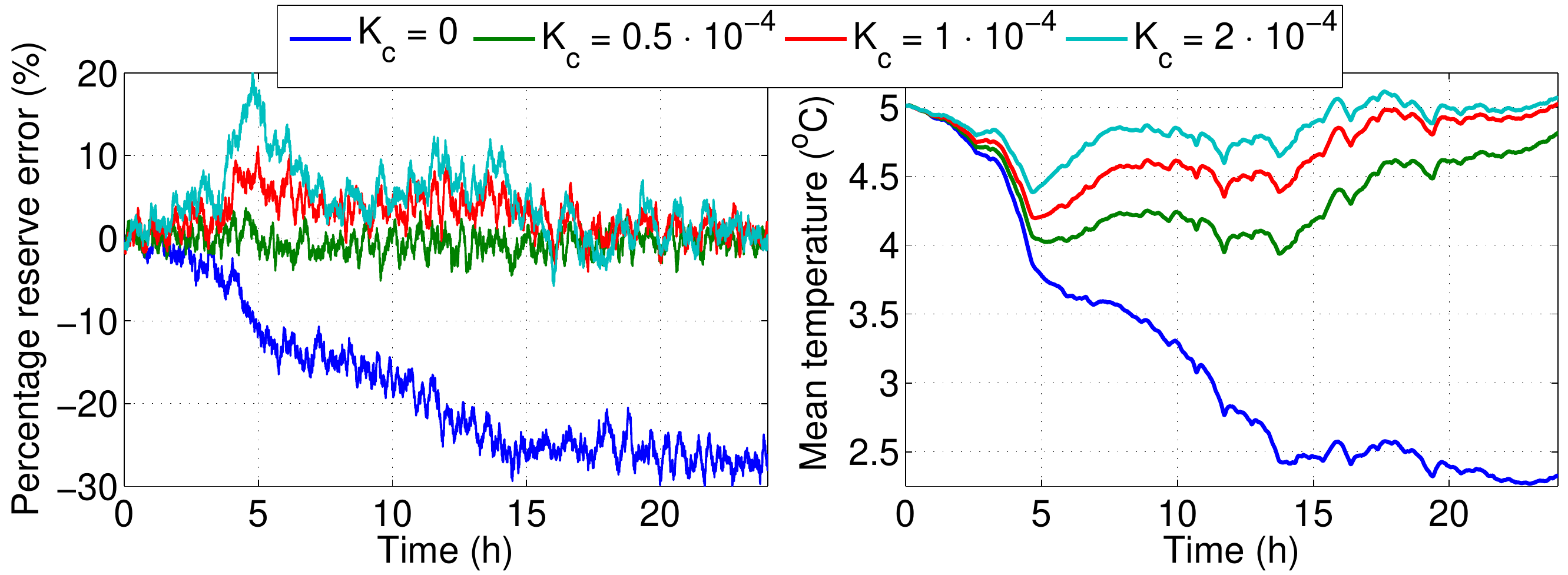}
\caption{Time-domain simulation results of reserve error and mean temperature deviation for four different values of the corrective temperature gain.} \label{fig:corrective_time_simulation}
\end{figure}

In this section we present results for the corrective temperature control loop. Figure~\ref{fig:corrective_time_simulation} compares the reserve error $e_{\text{r},t} = 100\cdot(P_{\text{d},t}-P_{\text{agg},t})/P_{\text{res}}$ and the temperature deviations for three $K_\text{c}$ values, as well as for a case without corrective control. The results correspond to simulations with a frequency deviation signal that is positively biased for the first $15$ hours (with an average bias of $\delta=0.0192$~Hz) and zero-mean for the rest of the day. If $K_\text{c}=0$, there is a steady-state error in temperature, which has a serious impact on user utility. At the same time, the temperature error affects the baseline power of the population, which results in a large PFC error and reduces the population's ability to provide PFC. On the other hand, all three non-zero gains bring the mean temperature close to the nominal value of $5^\circ\textrm{C}$ and eliminate the steady-state reserve error. As expected, higher $K_\text{c}$ values allow a faster temperature recovery but introduce oscillations in the baseline power, which translate into increased reserve errors at the beginning of the recovery period (hour $5$ in Fig.~\ref{fig:corrective_time_simulation}).

\begin{figure}[t]
\centering
\includegraphics[width=\textwidth]{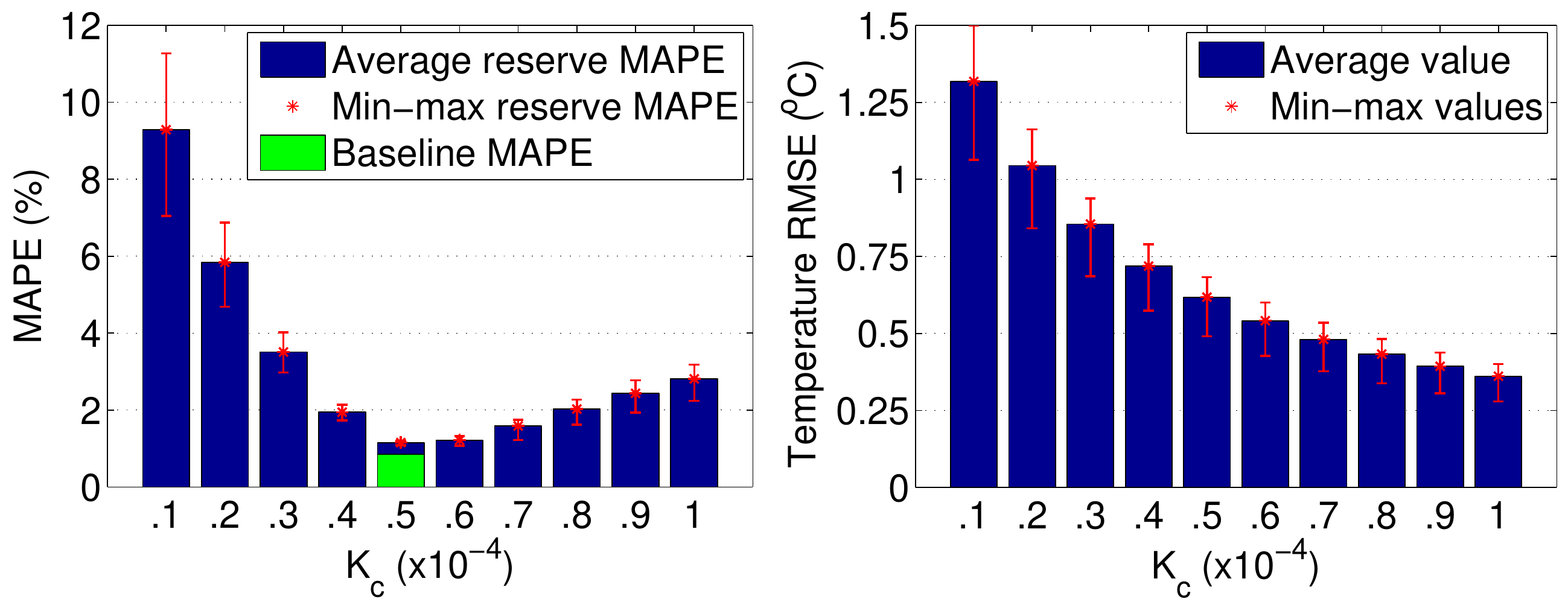}
\caption{Reserve and temperature errors for $10$ corrective temperature gains (averaged over $10$ simulations cases).} \label{fig:corrective_bar_plots}
\end{figure}

Results from $10$ day-long simulations for $K_\text{c}$ values in the range $[0.1,1] \cdot 10^{-4}$ are presented in Fig.~\ref{fig:corrective_bar_plots}. The frequency deviation is positively biased in $4$ of the days, negatively biased in another $4$ of the days, and zero-mean in the remaining $2$ days. The blue bars show the average values of reserve \ac{MAPE} and temperature \ac{RMSE}, whereas the red lines indicate the range spanned from the minimum to the maximum value. As expected, the reserve \ac{MAPE} is a convex function of $K_\text{c}$. The gain $K_\text{c}=0.5 \cdot 10^{-4}$ results in the minimum reserve \ac{MAPE} of $1.15\%$, which is very close to the baseline \ac{MAPE} of $0.85\%$ due to the system's natural dynamics (shown with the green bar in Fig.~\ref{fig:corrective_bar_plots}).

The results of Fig.~\ref{fig:corrective_time_simulation} are in agreement with Proposition~\ref{mean_temp_proposition}. With $\delta=0.0192$ Hz, $N_\text{ev}=15$~h, $N_\text{rec}=9$~h, $\varepsilon=1^\circ\textrm{C}$, $\epsilon=0.2^\circ\textrm{C}$, $D^\text{r}=0.15$, $\bar{P}_\text{n}=80$~W and $\bar{\beta}=4.4\cdot10^{-5}$, the lower bound on $K_\text{c}$ is $0.4863\cdot10^{-4}$. From Fig.~\ref{fig:corrective_time_simulation} we can see that $K_\text{c}=0.5\cdot10^{-4}$ is consistent with the design criteria: (i) the maximum temperature deviation is approximately $\varepsilon=1^\circ\textrm{C}$; and (ii) $N_\text{rec}=9$~h after the frequency bias has disappeared, the mean temperature is less than $\epsilon=0.2^\circ\textrm{C}$ away from the nominal value of $5^\circ\textrm{C}$. In addition, the results of Fig.~\ref{fig:corrective_bar_plots} are in agreement with \eqref{upper_Kc_bound}. For the considered values $\overline{\Delta{T}}=2^\circ\textrm{C}$, $\bar{T}_\text{a}=22^\circ\textrm{C}$, $\bar{\eta}\bar{R}\bar{P}_\text{n}=70^\circ\textrm{C}$ and $\Delta{t}=1$~s, \eqref{upper_Kc_bound} provides us with the upper bound $K_\text{c}=0.5004\cdot10^{-4}$. Therefore, the theoretical lower and upper bounds are very tight.

\subsection{Effect of Limited Thermostat Resolution} \label{thermostat_results}
\begin{figure}[t]
\centering
\includegraphics[width=\textwidth]{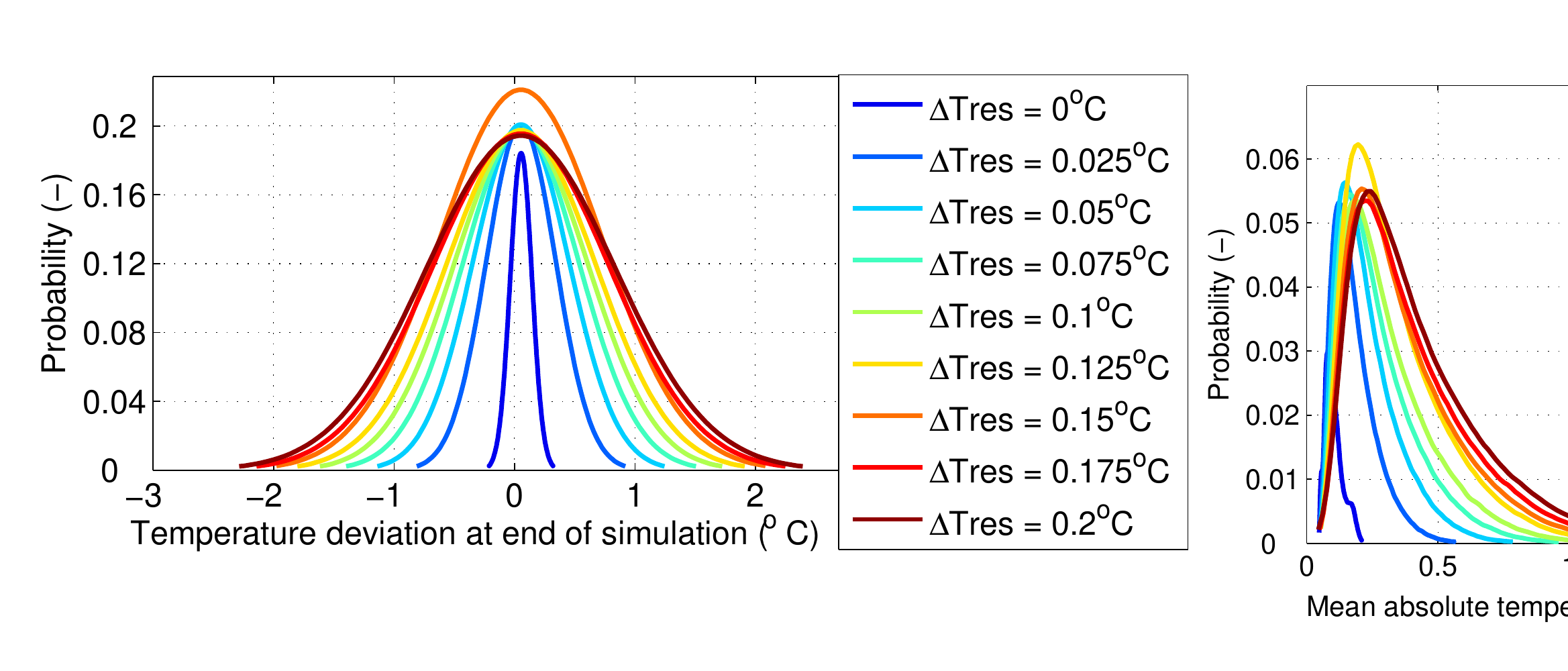}
\caption{The PDF of temperature deviations at the end of simulation.} \label{fig:thermostat_resolution_pdf}
\vspace{-0.2cm}
\end{figure}

To investigate the effect of minimum thermostat resolution $\Delta{T}_\text{res}$ on temperature deviations and control accuracy, we performed $10$ day-long simulations with $9$ different $\Delta{T}_\text{res}$ values and present the results in Fig.~\ref{fig:thermostat_resolution_pdf}. The plot shows the empirical \acp{PDF} of temperature deviations at the end of the simulation, which match closely normal distributions. The mean temperature deviation is independent of $\Delta{T}_\text{res}$, whereas the variance increases with $\Delta{T}_\text{res}$, i.e., larger deviations are possible for higher $\Delta{T}_\text{res}$ values. These simulation results are in agreement with Proposition~\ref{var_temp_proposition} and Remark~\ref{var_temp_remark}.

\begin{figure}[t]
\centering
\includegraphics[width=\textwidth]{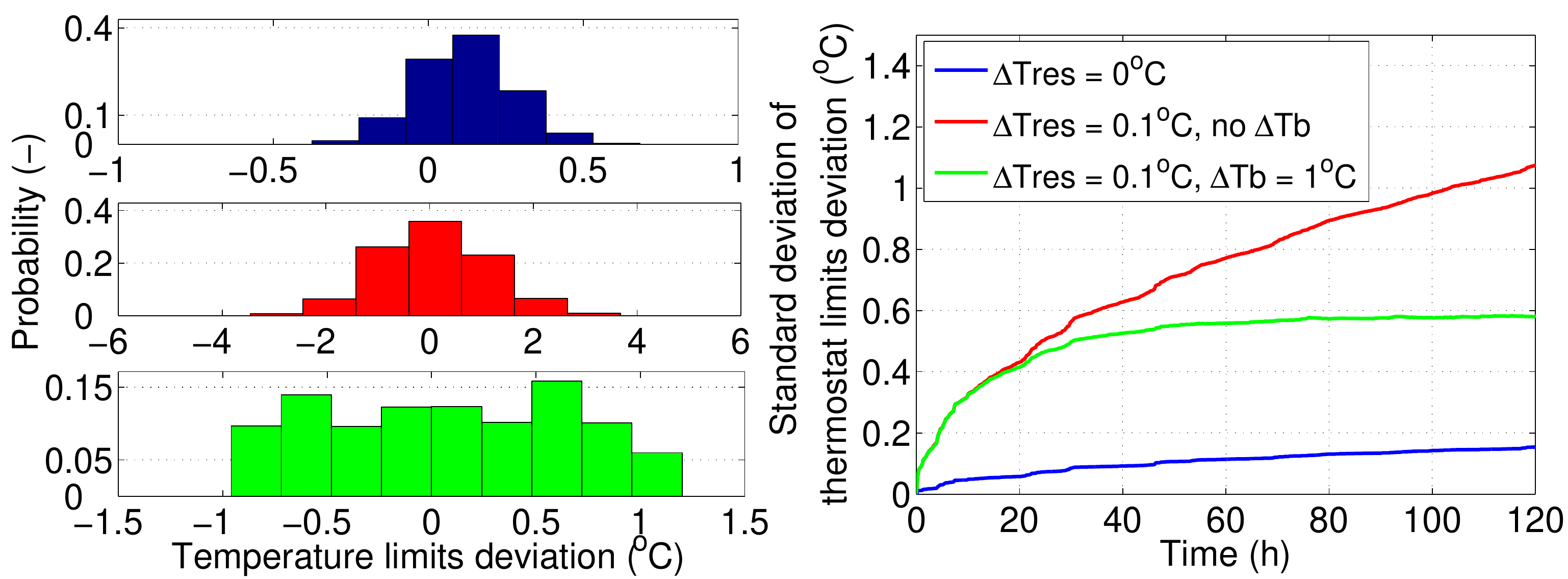}
\caption{The effect of limited thermostat resolution on the temperature deviations, with and without a temperature deviation bound.} \label{fig:thermostat_resolution_variance}
\end{figure}

Figure~\ref{fig:thermostat_resolution_variance} compares the following cases: (i) $\Delta{T}_\text{res}=0^\circ\textrm{C}$, (ii) $\Delta{T}_\text{res}=0.1^\circ\textrm{C}$ and no bound $\Delta{T}_\text{b}$, and (iii) $\Delta{T}_\text{b}=0.1^\circ\textrm{C}$ with $\Delta{T}_\text{b}=1^\circ\textrm{C}$. The left plot shows the histogram of temperature limit deviations at the end of a $5$-day simulation, whereas the right plot shows the evolution of the standard deviation of temperature limit changes over time. If $\Delta{T}_\text{res}=0^\circ\textrm{C}$ (blue color), the histogram resembles a normal distribution and the standard deviation slowly increases and settles to a small value, as explained in Remark~\ref{var_temp_remark}. If $\Delta{T}_\text{res}=0.1^\circ\textrm{C}$ and no $\Delta{T}_\text{b}$ is used (red color), the histogram resembles again a normal distribution but the standard deviation monotonically increases, as shown in Proposition~\ref{var_temp_proposition}. However, if a bound $\Delta{T}_\text{b}=1^\circ\textrm{C}$ is used (green color), the standard deviation is capped at a significantly lower value. Note that the bounding results in a more uniform PDF of temperature deviations in the population.

\begin{figure}[t]
\centering
\includegraphics[width=0.95\textwidth]{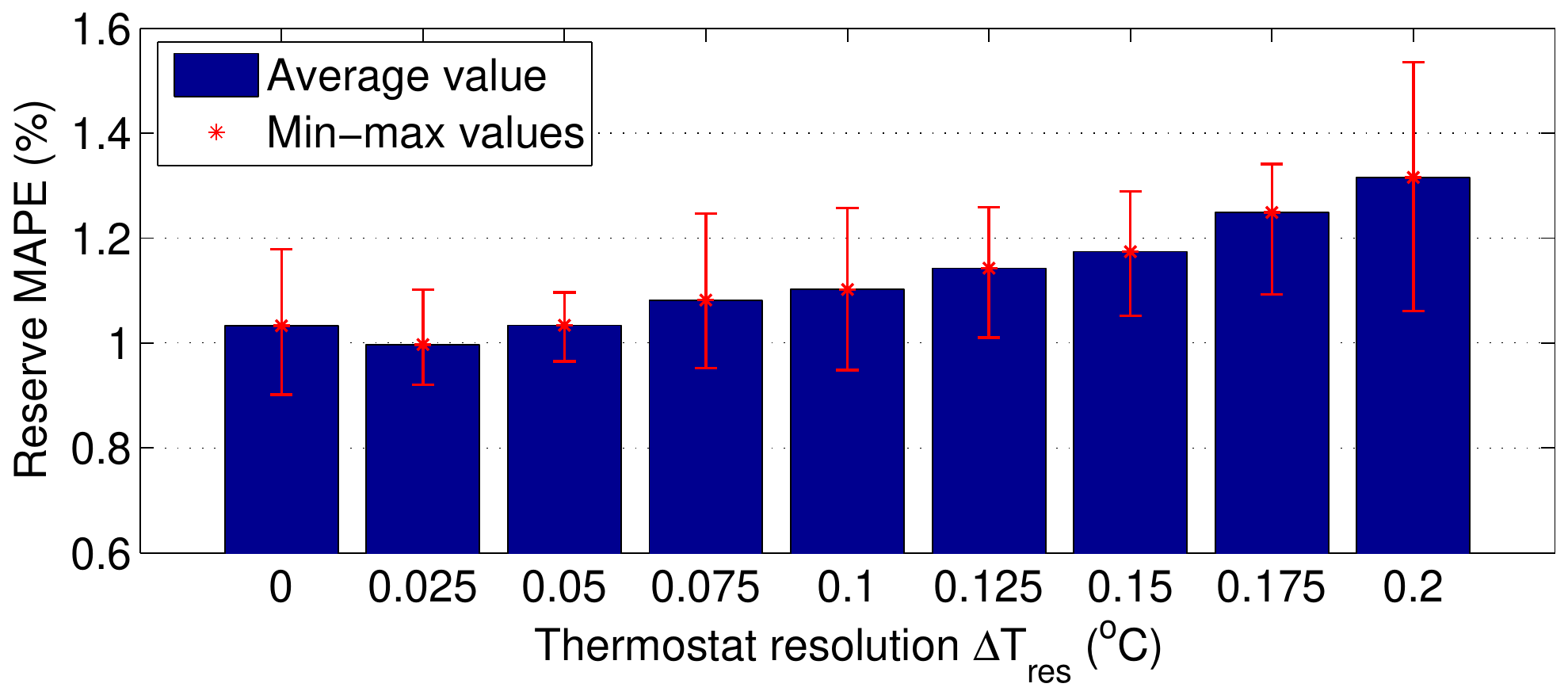}
\vspace{-0.4cm}
\caption{The dependence of reserve MAPE on the thermostat resolution (results averaged over $10$ simulations cases).}
\label{fig:thermostat_resolution_res_error}
\end{figure}

The effect of minimum thermostat resolution on control accuracy is shown in Fig.~\ref{fig:thermostat_resolution_res_error}. The blue bars are mean values of reserve MAPE over $10$ simulations, whereas the red lines indicate the minimum and maximum values. The reserve MAPE generally increases as $\Delta{T}_\text{res}$ increases but the increase is rather small, i.e., thermostat resolution limitations only slightly reduce control accuracy. In addition, we compared the reserve MAPE for $\Delta{T}_\text{res}=0.1^\circ\textrm{C}$, with and without a bound $\Delta{T}_\text{b}=1^\circ\textrm{C}$, and no observable reduction in control accuracy was found.

\subsection{Performance with Door Openings} \label{door_opening_results}
\begin{figure}[t]
\centering
\includegraphics[width=\textwidth]{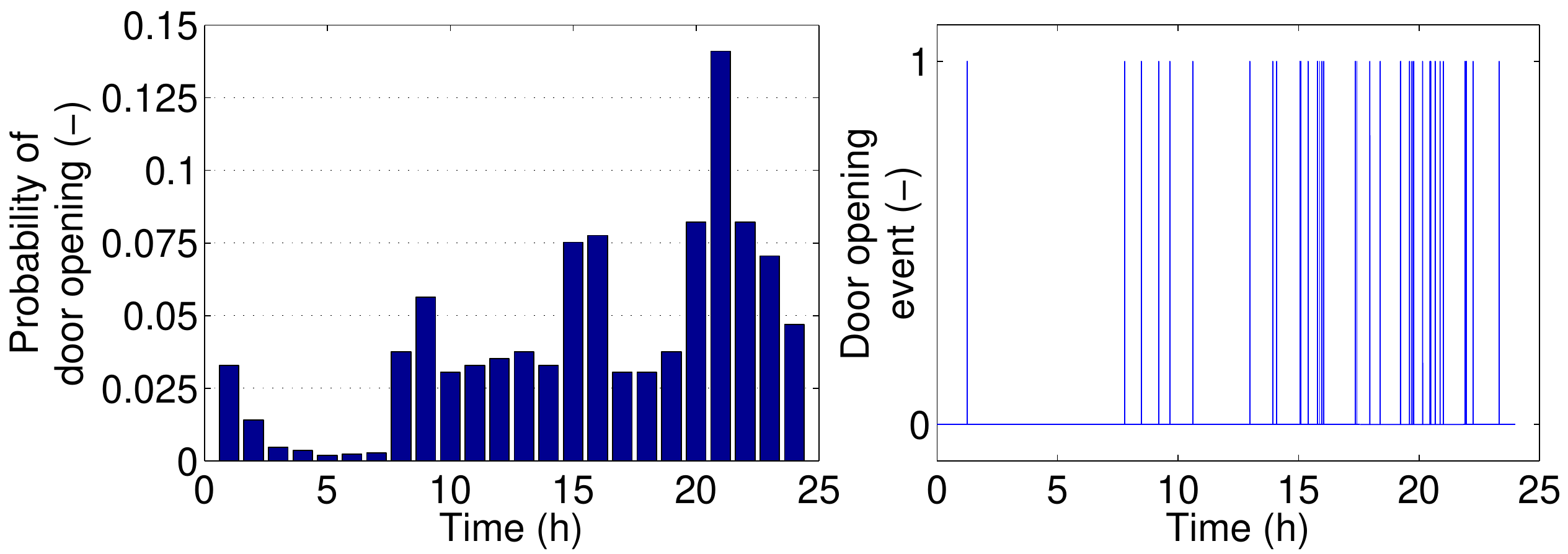}
\caption{Left: The histogram of door openings within a day. Right: Door opening events for a refrigerator generated based on the histogram.} \label{fig:door_opening_dist}
\end{figure}

In this section, we investigate the effect of refrigerator door openings on controller's performance. The left plot of Fig.~\ref{fig:door_opening_dist} shows the probability of door openings for each hour of the day based on data from \cite{parker1992openings}. The right plot shows a time series example of door opening events for a refrigerator, which is generated using the probability profile of the left plot.

\begin{figure}[t]
\centering
\includegraphics[width=\textwidth]{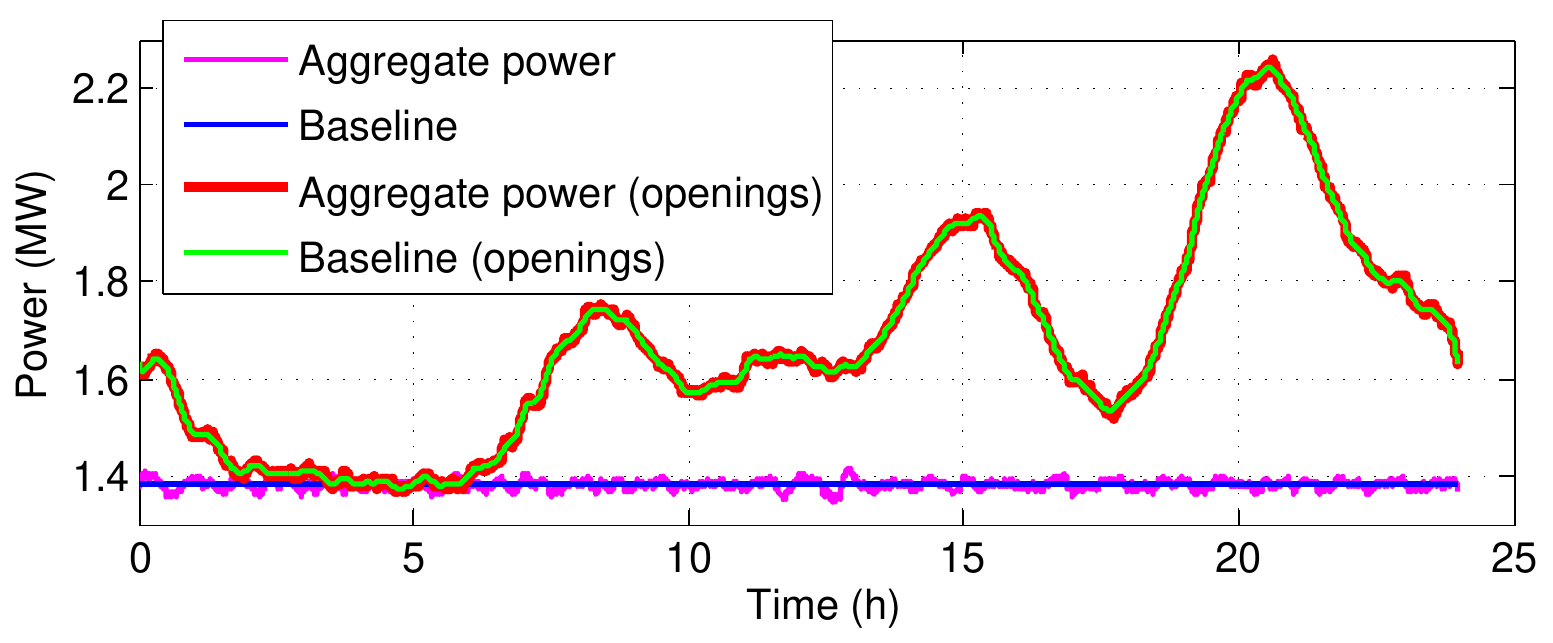}
\caption{$24$-hour baseline power with and without refrigerator door openings. The minimum baseline power of $1.38$~MW corresponds to a duty cycle of $0.25$, whereas the maximum power of $2.24$~MW to a duty cycle of $0.4$.} \label{fig:baseline_with_and_without_openings}
\end{figure}

We performed simulations for the base case scenario with and without door openings and compare the control performance in Table~\ref{table:MAPE with_openings}. Due to door openings, the reserve MAPE increases by $40\%$ for small frequency deviation bias and by $80\%$ for large frequency deviation bias. However, in both cases the error is still small and therefore the proposed controller is applicable despite the random door openings. The door openings smoothen the baseline power during load ramps, especially around $3$ pm and $7$ pm, and as a result the baseline MAPE decreases. To illustrate this, Fig.~\ref{fig:baseline_with_and_without_openings} shows the baseline power with and without door openings, as well as the aggregate power without PFC. Notice that there is less natural fluctuation around the baseline with door openings.

\begin{table}[t]
\renewcommand{\arraystretch}{1.05}
\caption{Baseline and Reserve MAPE with and without Door Openings}
\centering
\begin{tabular}{c c | c c | c c}
\hline
\multicolumn{2}{c|}{\begin{tabular}{@{}c@{}}Baseline \\ MAPE ($\%$)\end{tabular}} &
\multicolumn{2}{c|}{\begin{tabular}{@{}c@{}}Reserve  MAPE ($\%$) \\ Small bias\end{tabular}} & \multicolumn{2}{c}{\begin{tabular}{@{}c@{}}Reserve  MAPE ($\%$) \\ Large bias\end{tabular}}\\
\hline\
Without & With & Without & With &  Without & With \\
 $0.90$ & $0.47$ & $1.01$ & $1.41$ & $1.15$ & $2.06$ \\
\hline
\end{tabular}
\label{table:MAPE with_openings}
\end{table}

\begin{figure}[t]
\centering
\begin{minipage}{0.49\linewidth}
\centering \includegraphics[width=0.98\textwidth]{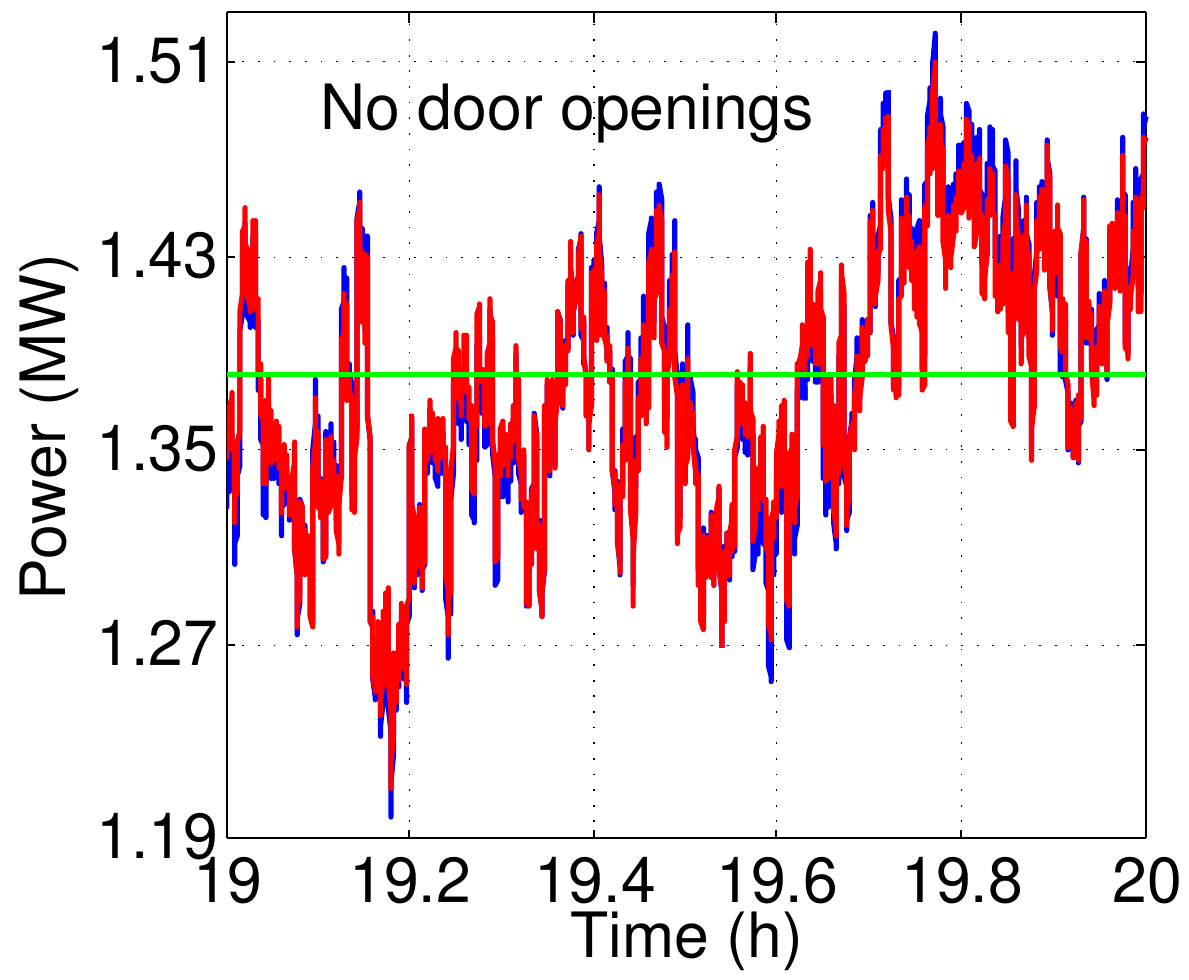}
\end{minipage}
\begin{minipage}{0.49\linewidth}
\centering \includegraphics[width=0.98\textwidth]{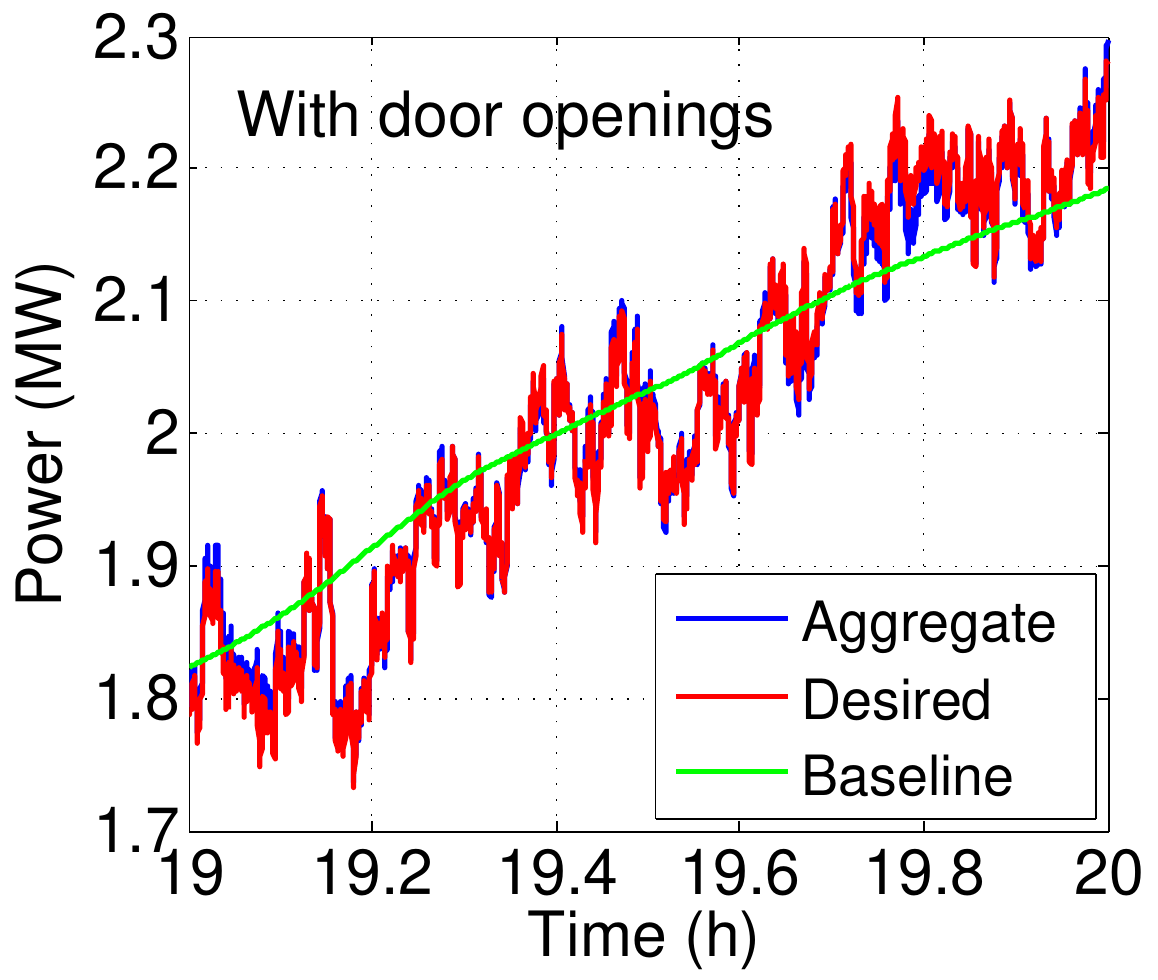}
\end{minipage}
\caption{PFC activation around the baseline without refrigerator door openings (left) and with refrigerator door openings (right) for a period of $1$ hour.} \label{fig:tracking_with_and_without_openings}
\end{figure}

Figure~\ref{fig:tracking_with_and_without_openings} shows the PFC activation (``desired power'') around the baseline,, with and without door openings, for a period of 1 hour. The baseline with door openings has an increasing trend, but the PFC activation is very close to that without door openings. Notice that the door openings do not have a significant effect on the tracking of the desired power. Therefore, if the door openings are considered in the control design as shown in the paper, then their effect on the tracking performance of the desired power will not be detrimental.

\section{Implementation Issues} \label{implementation_issues}
\subsection{Refrigerator Model}
Our analysis is based on the commonly used first-order freezer-less refrigerator model \eqref{eq:ode_off}, \eqref{eq:ode_on}. As shown in \cite{sossan2016grey}, using separate states for the temperature of the refrigerator compartment's air, refrigerator's content and evaporator, significantly improves model performance. Modeling a refrigerator with a freezer is possible following the approach of \cite{Keep2011switches}.

It is possible to integrate more detailed refrigerator models in the control design by modifying \eqref{on_time_est}, \eqref{off_time_est} and \eqref{duty_cycle_estimate1}. More importantly, the model of startup dynamics can be improved by considering the dependence on evaporator temperature when a refrigerator switches on. Since this temperature depends on the time elapsed since the previous on cycle, which is not constant when providing PFC, the profile of startup dynamics will likely be time-varying in reality.

\subsection{Noisy Frequency Measurements}
The authors of \cite{donnelly2010frequency} performed a correlation analysis of frequency measurements at the transmission level ($230$ kV) and at the wall outlet level ($120$ V). The results were promising for PFC and showed that the frequency measurements at the two voltage levels are highly correlated. Although it is important to investigate how measurement noise affects control performance, a moving average filter can be used to reduce the measurement noise in a practical implementation \cite{xu2011demand}.

\subsection{Cost and Economic Evaluation}
If $D^\text{r}=0.2$, our simulations show that $63,500$ refrigerators are necessary to provide $1$~MW of PFC reserves.\footnote{Assuming a refrigerator power rating of $80$~W and a duty cycle of $0.3$, the average power consumption of the aggregation of $63,500$ devices is $1.524$~MW. Therefore, the aggregation provides $100\cdot (1/1.524) \approx 65\%$ of the average power as reserve capacity.} Using a reserve capacity price of $21.5~\euro/\text{MW}/\text{h}$ (average price from January to September $2014$ in Switzerland), the total revenue from selling $1$~MW of capacity is $188,340~\euro/\text{year}$, which corresponds to a revenue per device as low as $3~\euro/\text{year}$. Since the average refrigerator lifetime is $14$ years, providing PFC reserves will be profitable only if the implementation cost is below $42~\euro$.

We expect the  micro-controller's cost to be a few $\euro$ due to its simplicity, whereas the frequency meter is the dominant cost driver of the proposed scheme. Reference \cite{xu2011demand} reports a cost of $20~\euro$ for a controller including the frequency meter, but without mentioning the meter's accuracy. According to our market search, a typical selling price in $2016$ for a frequency meter with accuracy $10$ mHz was approximately $40~\euro$, whereas the price explodes to approximately $600~\euro$ for an accuracy of $1$ mHz. Based on discussions with manufacturers of high-accuracy frequency meters, the cost is expected to drop to approximately $50~\euro$, if the meters are standardized and integrated in the refrigerators' electronic circuits.

Nevertheless, even with a reduced frequency meter cost, providing PFC with residential refrigerators will not be a very profitable investment due to the required large number of devices and the currently low reserve prices. In order to make customer participation in such programs more attractive, the refrigerator owners could enjoy electricity price rebates. In addition, performance-based reserve payments could be used in order to remunerate the refrigerators' higher accuracy and faster response in reserve provision. In fact, performance-based payments are already used for secondary frequency reserves in the PJM power market in the US.

The developed methods are also applicable to other cooling devices, for example, commercial refrigerators and freezers, refrigerated warehouses and air conditioners, with no or minimal modifications. Due to the much higher power rating of these devices, a significantly smaller aggregation will be needed for the same amount of reserves. Therefore, we expect the revenue for each participant to be notably higher, which makes this application case financially more attractive compared with residential refrigerators.

\section{Concluding Remarks} \label{conclusions}
In this paper, we presented a decentralized stochastic control scheme to enable PFC reserve provision from aggregations of refrigerators. The control is based on probabilistic switching to avoid load synchronization, and it accounts for refrigerator startup dynamics and lockout constraints. Reliable reserve provision during persistent frequency deviations is achieved by resetting the thermostat temperature limits. Furthermore, the controller is robust to biased frequency deviations, limited thermostat resolution and door openings. Extensive simulation results show that the control scheme allows a refrigerator aggregation to provide fast and reliable PFC without communication. In the future, we plan to investigate the controller's robustness to excessive compressor locking, and perform dynamic simulation studies in a two-area power system model.

\section*{Appendix: Proof of Proposition~\ref{startup_proposition}} \label{proof_propo1}
\begin{proof}
The estimated aggregate power at time step $t$ after reserve activation is
\begin{align} \label{est_Pagg}
\hat{P}_{\text{agg},t} = N_\text{r,a} \bar{P}_\text{n} \bar{S}_{\text{u}}(t)~,
\end{align}
where $N_\text{r,a}$ is the number of activated devices. Assuming independence of random variables $P_{\textrm{n},i}$, $u_i$ and $N_{\textrm{s},i}$, the expected actual aggregate power is computed as
\begin{align}
\bar{P}_{\text{agg},t} &= \Mean\Big[\sum\nolimits_{i=1}^{N_\text{r,a}} P_{i,t}\Big] = \sum\nolimits_{i=1}^{N_\text{r,a}}\Mean\left[P_{i,t}\right]\nonumber
\end{align}
\begin{align}
&= \sum_{i=1}^{N_\text{r,a}}\Mean\Bigg[P_{\textrm{n},i} \Bigg[1+ u_i \bigg[1-\frac{t}{N_{\textrm{s},i}}\bigg]_+\Bigg]\Bigg]\nonumber
\end{align}
\begin{align}\label{exp_Pagg}
&= \sum_{i=1}^{N_\text{r,a}} \bar{P}_\text{n}\Bigg[1+\bar{u}\cdot\Mean\bigg[1-\frac{t}{N_{\textrm{s},i}}\bigg]_+\Bigg]~.
\end{align}

\underline{\emph{Part 1:}} If $t \leq N_\text{s,min}$, then $[\cdot]_+$ is redundant and \eqref{exp_Pagg} gives
\begin{align}\label{exp_Pagg2}
&\bar{P}_{\text{agg},t} = \sum\nolimits_{i=1}^{N_\text{r,a}} \bar{P}_\text{n}\bigg[1+\bar{u}\bigg[1-\Mean\bigg[\frac{t}{N_{\textrm{s},i}}\bigg]\bigg]\bigg] \nonumber\\
&\leq \hspace{-0.06cm} \sum\nolimits_{i=1}^{N_\text{r,a}} \bar{P}_\text{n}\bigg[1+\bar{u}\bigg(1-\frac{t}{\bar{N}_\textrm{s}}\bigg)\bigg]\hspace{-0.06cm}=\hspace{-0.06cm} N_\text{r,a} \bar{P}_\text{n} \bar{S}_{\text{u}}(t) \hspace{-0.06cm}=\hspace{-0.06cm} \hat{P}_{\text{agg},t},
\end{align}
where Jensen's inequality $\phi(\Mean[N_{\textrm{s},i}]) \leq \Mean[\phi(N_{\textrm{s},i})]$ is invoked with $\phi(N_{\textrm{s},i})=t/N_{\textrm{s},i}$, which is a convex function since $t, N_{\textrm{s},i} \geq 0$. Thus, $\bar{P}_{\text{agg},t} \leq \hat{P}_{\text{agg},t}$ holds for $t \leq N_\text{s,min}$.

\underline{\emph{Part 2:}} For $t \geq \bar{N}_\textrm{s}$, we obtain $\hat{P}_{\text{agg},t}=0$ from \eqref{Su_peak_dyn} and \eqref{est_Pagg}. Since the expected aggregate power $\bar{P}_{\text{agg},t}$ is non-negative, $\hat{P}_{\text{agg},t}$ is clearly a lower bound.

If $N_{\textrm{s},i}$ follows the uniform PDF $N_{\textrm{s},i}$$\scriptstyle\sim$$[N_\text{s,min}, N_\text{s,max}]$, then for a given $t \in (N_\text{s,min}, \bar{N}_\textrm{s})$ the probability that a randomly selected refrigerator contributes to the aggregate startup dynamics is equal to $\Prob\left(N_{\textrm{s},i} \geq t\right) = \left(N_\text{s,max}-t\right)/\left(N_\text{s,max}-N_\text{s,min}\right)$. Let us denote by $\tilde{N}_{\textrm{s},i}$ the random startup duration that follows the uniform PDF $\tilde{N}_{\textrm{s},i}$$\scriptstyle\sim$$[t, N_\text{s,max}]$. The expectation from the last step of \eqref{exp_Pagg} can be written as
\begin{subequations}
\begin{align}
&\Mean\bigg[1-\frac{t}{N_{\textrm{s},i}}\bigg]_+ = \Mean\bigg[1-\frac{t}{\tilde{N}_{\textrm{s},i}}\bigg] \cdot \frac{N_\text{s,max}-t}{N_\text{s,max}-N_\text{s,min}} \\
&\leq \bigg(1-\frac{t}{\Mean[ \tilde{N}_{\textrm{s},i} ]}\bigg) \cdot \frac{N_\text{s,max}-t}{N_\text{s,max}-N_\text{s,min}} \label{jensens_use}\\
&= \bigg[1-\frac{t}{0.5\cdot(t+N_\text{s,max})}\bigg] \cdot \frac{N_\text{s,max}-t}{N_\text{s,max}-N_\text{s,min}}\label{uniform_use}\\
&= \frac{\left(N_\text{s,max}-t\right)^2}{\left(N_\text{s,max}-N_\text{s,min}\right)\cdot\left(N_\text{s,max}+t\right)}~.
\end{align}
\end{subequations}
We used Jensen's inequality in \eqref{jensens_use} and the fact that $\tilde{N}_{\textrm{s},i}$ follows the uniform distribution $\tilde{N}_{\textrm{s},i}$$\scriptstyle\sim$$[t, N_\text{s,max}]$ in step \eqref{uniform_use}.

The inequality $\hat{P}_{\text{agg},t} \geq \bar{P}_{\text{agg},t}$ holds if and only if
\begin{align}
1-\frac{t}{\bar{N}_\textrm{s}} &\geq \Mean\bigg[1-\frac{t}{N_{\textrm{s},i}}\bigg]_+ \Leftrightarrow \\
\frac{N_\text{s,max}+N_\text{s,min}-2t}{N_\text{s,max}+N_\text{s,min}} &\geq \frac{\left(N_\text{s,max}-t\right)^2}{\left(N_\text{s,max}\hspace{-0.05cm}-\hspace{-0.05cm}N_\text{s,min}\right)\hspace{-0.05cm}\cdot\hspace{-0.06cm}\left(N_\text{s,max}\hspace{-0.05cm}+\hspace{-0.05cm}t\right)}~,\label{complicated_inequality}
\end{align}
where we used the fact that $\bar{N}_\textrm{s}=0.5\cdot(N_\text{s,max}+N_\text{s,min})$. With some effort \eqref{complicated_inequality} can be equivalently written as
\begin{align}
(N_\text{s,min}&-3N_\text{s,max}) \cdot t^2 + (N_\text{s,max}^2 + 4N_\text{s,min}N_\text{s,max}-N_\text{s,min}^2) \cdot t \nonumber\\
&- (N_\text{s,min}^2N_\text{s,max}+N_\text{s,min}N_\text{s,max}^2) \geq 0~.\label{quadratic_form}
\end{align}
The quadratic form in \eqref{quadratic_form} has two real solutions $t_1$ and $t_\text{lim}$. The first one is $t_1=N_\text{s,min}$ by inspection. For the second one, we use Vieta's formula $t_1 \cdot t_\text{lim} = c/a$ with $a=N_\text{s,min}-3N_\text{s,max}$ and $c=- (N_\text{s,min}^2N_\text{s,max}+N_\text{s,min}N_\text{s,max}^2)$ and find
\begin{align}
t_\text{lim} = \frac{N_\text{s,max} \left(N_\text{s,min}+N_\text{s,max}\right)}{3N_\text{s,max}-N_\text{s,min}}~.
\end{align}
The quadratic inequality \eqref{quadratic_form} holds for $t\in [t_1, t_\text{lim}]$. Therefore, if $t > N_\text{s,min}$ the inequality $\hat{P}_{\text{agg},t} \geq \bar{P}_{\text{agg},t}$ holds for $t\leq t_\text{lim}$.
\end{proof}

\section*{Appendix: Proof of Proposition~\ref{mean_temp_proposition}} \label{proof_propo3}
\begin{proof}
\underline{\emph{Part 1:}} For $0 \leq t \leq N_\text{ev}+1$ using $\bar{T}_0 = \bar{T}_{\text{nom}}$ and \eqref{Df_signal}, equation \eqref{mean_temp_est_simple_rec} gives
\begin{align}
\bar{T}_t &= \bar{T}_{\text{nom}} - \gamma \delta \cdot \sum\nolimits_{k=0}^{t-1} \lambda^k
= \bar{T}_{\text{nom}} - \gamma \delta \cdot \frac{1-\lambda^{t-1}}{1-\lambda}~, \label{geom_series}
\end{align}
where a standard property of geometric series is invoked. The minimum of \eqref{geom_series} is obtained for $t=N_\text{ev}+1$ and thus if $\bar{T}_{\text{nom}}-\varepsilon \leq \bar{T}_t \leq \bar{T}_{\text{nom}}+\varepsilon$ holds for $t=N_\text{ev}+1$ then it holds $\forall t$. Assuming $\delta>0$ without loss of generality, it should hold
\begin{align}\label{lambda_ineq1_proof}
\bar{T}_{t=N_\text{ev}+1} \geq \bar{T}_{\text{nom}}-\varepsilon ~~\Rightarrow~~ -\gamma \delta \cdot \frac{1-\lambda^{N_\text{ev}}}{1-\lambda} \geq -\varepsilon~.
\end{align}

\underline{\emph{Part 2:}} For $t \geq N_\text{ev}+1$, $\Delta{f}_t=0$ and the mean temperature recovers towards $\bar{T}_{\text{nom}}$. We reset the time index such that $t^\prime=0 \Leftrightarrow t=N_\text{ev}+1$. The new initial condition is $\bar{T}_{(t^\prime=0)} = \bar{T}_{(t=N_\text{ev}+1)}$ and is computed from \eqref{geom_series}. Equation \eqref{mean_temp_est_simple_rec} gives
\begin{align}
\bar{T}_{t^\prime} &= \bar{T}_{\text{nom}} - \gamma \delta \cdot \frac{1-\lambda^{N_\text{ev}}}{1-\lambda} \cdot \lambda^{t^\prime}~. \label{geom_series2}
\end{align}

At time $t^\prime=N_\text{rec}$, $\bar{T}_{\text{nom}}-\epsilon \leq \bar{T}_{N_\text{rec}} \leq \bar{T}_{\text{nom}}+\epsilon$ must hold. Assuming $\delta>0$ without loss of generality, we get
\begin{align}\label{lambda_ineq2_proof}
\bar{T}_{t^\prime=N_\text{rec}} \geq \bar{T}_{\text{nom}}-\epsilon \Rightarrow -\gamma \delta \cdot \lambda^{N_\text{rec}} \cdot \frac{1-\lambda^{N_\text{ev}}}{1-\lambda} \geq -\epsilon~.
\end{align}
Equations \eqref{lambda_ineq1_proof} and \eqref{lambda_ineq2_proof} are equivalent to \eqref{lambda_ineq1} and \eqref{lambda_ineq2}.
\end{proof}

\section*{Appendix: Proof of Proposition~\ref{var_temp_proposition}} \label{proof_propo4}
\begin{proof}
The change in temperature limits of a randomly selected refrigerator $i$ can be expressed as
\begin{align}
\Delta{T}_{\text{lim},i,t} = \sum\nolimits_{k=0}^{t-1} \Delta{T}_\text{res} \cdot \zeta_k \cdot z_k ~,
\end{align}
where $\zeta_k$ is equal to $1$ or $-1$ depending on the sign of $\Delta{T}_{\text{lim},k}$, and $z_k$ is the Bernoulli random variable
\begin{align} \label{Bernoulli_var_temp}
f_{\text{z}}(z_k) = \begin{cases}
                        \varrho_k, & \textrm{if} \hspace{2mm} z_k=1 \vspace{-0.15cm}\\			
                        1-\varrho_k, & \textrm{if} \hspace{2mm} z_k=0
                    \end{cases}~.
\end{align}
Using the fact $\Mean\left[z_k\right]=\varrho_k$, the mean value is computed as
\begin{align}
\Mean\left[\Delta{T}_{\text{lim},i,t}\right] &= \sum\nolimits_{k=0}^{t-1}\Delta{T}_\text{res} \zeta_k \Mean\left[z_k\right] = \sum\nolimits_{k=0}^{t-1} \Delta{T}_\text{res} \zeta_k \varrho_k \nonumber \\
&= \sum\nolimits_{k=0}^{t-1} \zeta_k \left|\Delta{T}_{\text{lim},k}\right| = \sum\nolimits_{k=0}^{t-1} \Delta{T}_{\text{lim},k}~.\label{mean_comp}
\end{align}

The variance is computed as
\begin{subequations}
\begin{align}
&\mathrm{Var}\left[\Delta{T}_{\text{lim},i,t}\right] = \mathrm{Var}\Big[\sum\nolimits_{k=0}^{t-1} \Delta{T}_\text{res} \cdot \zeta_k \cdot z_k\Big] \label{var_comp1}\\
&= \sum_{k=0}^{t-1} \Delta{T}_\text{res}^2 \hspace{-0.02cm}\cdot\hspace{-0.02cm} \zeta_k^2 \hspace{-0.02cm}\cdot\hspace{-0.02cm} \mathrm{Var}\left[z_k\right]
= \sum_{k=0}^{t-1} \Delta{T}_\text{res}^2 \hspace{-0.02cm}\cdot\hspace{-0.02cm} \varrho_k \hspace{-0.02cm}\cdot\hspace{-0.02cm} (1-\varrho_k) \label{var_comp2}\\
&= \sum\nolimits_{k=0}^{t-1} \Delta{T}_\text{res}^2 \cdot \frac{|\Delta{T}_{\text{lim},k}|}{\Delta{T}_\text{res}} \cdot \Big(1-\frac{|\Delta{T}_{\text{lim},k}|}{\Delta{T}_\text{res}}\Big)\label{var_comp3}\\
&= \sum\nolimits_{k=0}^{t-1} |\Delta{T}_{\text{lim},k}| \cdot \big(\Delta{T}_\text{res}-|\Delta{T}_{\text{lim},k}|\big) \label{var_comp4}~,
\end{align}
\end{subequations}
where the Bernoulli random variables $z_k$ and $z_j$ are assumed to be uncorrelated $\forall k\hspace{-0.02cm}\neq\hspace{-0.02cm} j$ in \eqref{var_comp1}, $\mathrm{Var}\left[z_k\right]\hspace{-0.02cm}=\hspace{-0.02cm}\varrho_k \hspace{-0.02cm}\cdot\hspace{-0.02cm} (1-\varrho_k)$ is used in \eqref{var_comp2}, and the definition of $\varrho_k$ from \eqref{prob_temp_resetting} is used in \eqref{var_comp3}.
\end{proof}

\section*{Appendix: Proof of Proposition~\ref{door_opening_prop}} \label{proof_propo5}
\begin{proof}
Denote by $E^\text{d}_{\text{op}}$ and $E^\text{d}_{\text{cl}}$ the refrigerator's daily energy consumption with and without door openings, respectively. Consider a refrigerator subject to door openings and denote by $\bar{T}_\text{op}$ and $\bar{T}_\text{cl}$ the time-averaged temperature while the door is open and closed, respectively. The refrigerator operates on average $\mu_\text{op} \cdot \mu_\text{d}$ seconds with open door, and $N^\text{d}-\mu_\text{op} \cdot \mu_\text{d}$ seconds with closed door. The refrigerator's time constant with closed door is $\alpha=1/(RC)$, whereas with open door it is $\alpha_{\text{op}}=1/(R_{\text{op}}C)$.

The refrigerator's energy consumption is equal to the integral of instantaneous thermal losses over time because its temperature is regulated. Assuming that $\bar{T}_\text{cl}$ and $\bar{T}_\text{op}$ are constant, the daily energy consumption with and without door openings can be approximated using \eqref{eq:ode_off}, \eqref{eq:ode_on} as
\begin{align}
&E^\text{d}_{\text{cl}} \approx (\alpha/\beta) \cdot (T_\text{a}-\bar{T}_\text{cl}) \cdot N^\text{d} \label{energy_closed}\\
&E^\text{d}_{\text{op}} \approx (\alpha/\beta) \cdot (T_\text{a}-\bar{T}_\text{cl}) \cdot (N^\text{d}-\mu_\text{op} \mu_\text{d}) +\\
&\quad\quad\quad\quad\quad (\alpha_{\text{op}}/\beta) \cdot (T_\text{a}-\bar{T}_\text{op}) \cdot \mu_\text{op} \mu_\text{d}, \label{energy_open}
\end{align}
whereas $E^\text{d}_{\text{op}}$ and $E^\text{d}_{\text{cl}}$ are related according to
\begin{align} \label{energy_closed_and_open}
E^\text{d}_{\text{op}} = (1+\xi) \cdot E^\text{d}_{\text{cl}}~.
\end{align}

Combining \eqref{energy_closed}, \eqref{energy_open} and \eqref{energy_closed_and_open}, we get
\begin{align} \label{energy_aux}
&(\alpha/\beta) (T_\text{a}-\bar{T}_\text{cl}) (N^\text{d}-\mu_\text{op} \mu_\text{d}) +  (\alpha_{\text{op}}/\beta) (T_\text{a}-\bar{T}_\text{op}) \mu_\text{op} \mu_\text{d} = \nonumber\\
&\quad\quad\quad\quad\quad\quad\quad(1+\xi) \cdot (\alpha/\beta) \cdot (T_\text{a}-\bar{T}_\text{cl}) \cdot N^\text{d}~,
\end{align}
which when solved for $\alpha_{\text{op}}$ gives
\begin{align}
\alpha_{\text{op}} &= \alpha \cdot \frac{T_\text{a}-\bar{T}_\text{cl}}{T_\text{a}-\bar{T}_\text{op}} \cdot \bigg(1+\frac{N^\text{d}}{\mu_\text{op} \cdot \mu_\text{d}}\cdot\xi\bigg) ~~\Rightarrow \label{aopen1}\\
\alpha_{\text{op}} &\geq \alpha \cdot [1+(N^\text{d}\cdot\xi)/(\mu_\text{op}\cdot\mu_\text{d})]~,\label{aopen2}
%\alpha_{\text{op}} &\geq \alpha \cdot \bigg(1+\frac{N^\text{d}}{\mu_\text{op} \cdot \mu_\text{d}}\cdot\xi\bigg)~,\label{aopen2}
\end{align}
because $\bar{T}_\text{cl}\leq \bar{T}_\text{op}$. Inequality \eqref{Ropen} now directly follows from \eqref{aopen2} using $\alpha=(1/RC)$ and $\alpha_{\text{op}}=(1/R_{\text{op}}C)$.
\end{proof}

\bibliographystyle{IEEEtran}
\bibliography{biblio_EV}

\end{document}